\documentclass[a4paper,11pt]{elsarticle}
\usepackage{amsfonts}
\usepackage{amssymb}
\usepackage{amsmath}

\usepackage[utf8]{inputenc}
\usepackage[french]{babel}
\usepackage[T1]{fontenc}

\usepackage[utf8]{inputenc}
\usepackage{amsfonts}
\usepackage{amsthm}
\usepackage{mathtools}
\usepackage{enumitem}
\usepackage{esint}
\usepackage{dsfont}

\theoremstyle{plain}
\newtheorem{thm}{Theorem}[section]
\theoremstyle{definition}
\newtheorem{defn}[thm]{Definition}
\theoremstyle{plain}
\newtheorem{lem}[thm]{Lemma}
\theoremstyle{plain}

\theoremstyle{plain}

\theoremstyle{remark}
\newtheorem*{rem*}{Remark}

\usepackage{geometry}
\begin{document}

\begin{frontmatter}

  \title{Construction of Banach frames and atomic decompositions of anisotropic Besov spaces}
  
  \author{Dimitri Bytchenkoff $^{1, 2, *}$}
  
  \address{$^{1}$Technische Universit\"at Berlin, Institut f\"ur Mathematik, Stra{\ss}e des 17. Juni 136, 10623 Berlin, Deutschland}
    \address{$^{2}$Universit{\'e} de Lorraine, Laboratoire d'Energ{\'e}tique et de M{\'e}canique Th{\'e}orique et Appliqu{\'e}e, 2 avenue de la For{\^e}t de Haye, 54505 Vandoeuvre-l{\`e}s-Nancy, France}
   
   \begin{abstract}
We construct generalised shift-invariant systems of functions of several real variables for anisotropic Besov spaces that can be generated by the decomposition method using any given expansive matrix and establish the conditions on those systems under which they will constitute Banach frames or sets of atoms for the anisotropic homo- or inhomogeneous Besov spaces.

   \end{abstract}
   

\end{frontmatter}

\let\thefootnote\relax\footnotetext{* Corresponding author; e-mail address: dimitri.bytchenkoff@univ-lorraine.fr}

\section{Introduction}

\noindent Besov spaces, originally constructed by the approximation method \cite{Besov_1961}, play an extremely important role in the theory of differentiable functions of several real variables as they, on the one hand, constitute a closed system with respect to embedding theorems and are, on the other hand, closely related to Sobolev spaces \cite{Sobolev_1950}. Along with Sobolev, Besov spaces are an integral part of the embedding theory, which studies connexions between differential properties of functions in different metrics \cite{Sobolev_1950, Nikolski_1969, Besov_1975}. Harmonic analysis uses either bases or  frames and sets of atoms to decompose functions of a function space into basic building blocks or synthesise them from those blocks. The isotropic Besov spaces are known \cite{Triebel_2006} to have orthonormal bases made of wavelets \cite{Daubechies_1992}. Herein we shall use the innovative approach reported in \cite{Voigtlaender_2016, Voigtlaender_2017} to construct Banach frames and sets of atoms for the anisotropic homo- and inhomogeneous Besov spaces introduced in \cite{Bownik_2005} as decomposition spaces \cite{Feichtinger_1985}. Similar results, although formulated and achieved rather differently, were reported in \cite{Bownik_2005} and \cite{Borup_2007}. The work \cite{Bownik_2005} was, in its turn, an ingenious generalisation of the ideas developed in \cite{Frazier_1985} and \cite{Frazier_1989}.

This report is structured as follows. The elements of the decomposition method  essential to the present work are outlined at the beginning of Section \ref{Construction of anisotropic Besov spaces by the decomposition method}. The definitions of the anisotropic homo- and inhomogeneous Besov spaces viewed as decompositions spaces follow in Subsections \ref{Construction of anisotropic homogeneous Besov spaces} and
\ref{Construction of anisotropic heterogeneous Besov spaces}. The notion of the Banach frame and that of the set of atoms for the decomposition space are reminded at the beginning of Section \ref{Construction of Banach frames for and atomic decompositions of the decomposition spaces}. In the same section we give the statements of the two theorems that we shall use to construct Banach frames and atomic decompositions of anisotropic Besov spaces. In Subsection \ref{Banach frames for and atomic decompositions of homogeneous Besov spaces} of Section \ref{Construction of Banach frames for and atomic decompositions of the anisotropic Besov spaces} the set of anisotropic homogeneous Besov wavelets is defined as a generalised shift-invariant system and the conditions are established under which this set will form a Banach frame or a set of atoms for the anisotropic homogeneous Besov space. Finally anisotropic inhomogeneous wavelets are defined in Subsection \ref{Banach frames for and atomic decompositions of heterogeneous Besov spaces}, along with the conditions under which they will be a Banach frame or a set of atoms for the anisotropic inhomogeneous Besov space.

Throughout this report the information is provided as it is first needed. All definitions and propositions borrowed from other works are supplied with corresponding references. All the new propositions are followed by their proofs.

\section{Construction of anisotropic Besov spaces by the decomposition method} \label{Construction of anisotropic Besov spaces by the decomposition method}

\noindent In this study we shall concern ourselves with the anisotropic Besov spaces as they were defined in \cite{Bownik_2005} by the decomposition method, initially reported in \cite{Feichtinger_1985} and further developed in \cite{Borup_2007}. This definition involves three basic building blocks, namely an almost structured cover $\mathcal{Q}$ of the open subset of the frequency space,
a regular partition of unity on the subset subordinate to the covering $\mathcal{Q}$ and a $\mathcal{Q}$-moderate weight. Here are the definitions of these three notions and those of the decomposition space and its reservoir.

\begin{defn}\label{almost_struct_cover}
The set $\mathcal{Q} = \{ Q_i \}_{i \in I}$ is called an $almost$  $structured$ $cover$ of the open subset $O$ of $\mathbb{R}^d$ where $d \in \mathbb{N}$, if
\begin{enumerate}
  \item $\mathcal{Q}$ is $admissible$, i.e. the number of elements in the sets $\left\{i' \in I : Q_{i'} \cap Q_i \neq \emptyset \right\}$ is uniformly bounded for all $i \in I$;
    
  \item there is a set $\{ T_i \bullet + b_i \}_{i \in I}$ of invertible affine-linear maps and finite sets $\{ Q'_n \}_{n=1}^N$ and $\{P_n\}_{n=1}^N$ of non-empty open and bounded subsets $Q'_i$ and $P_i$ of $\mathbb{R}^d$ such that
   \begin{enumerate}                 
      \item $\overline{P_n} \subset Q'_n$ for all $1 \leqslant n \leqslant N$;
      \item for each $i \in I$ there is such an $n_i \in \{ 1, ..., N \}$ that $Q_i = T_i \, Q'_{n_i} + b_i$;        
      \item there is such a constant $C > 0$ that $\|T_i^{-1}T_{i'}\| \leqslant C$ for all such
          $i$ and $i' \in I$ that $Q_i \cap Q_{i'} \neq \emptyset$; and
      \item $O \subset \bigcup_{i \in I} (T_i P_{n_i} + b_i)$.
    \end{enumerate} 
\end{enumerate}
\end{defn}

\begin{defn}\label{reg_partition_unity}
Let $\mathcal{Q} = \{ Q_i \}_{i \in I}$ and $\{ T_i \bullet + b_i \}_{i \in I}$  be an almost structured cover of the open subset $O$ of $\mathbb{R}^d$ and the set of invertible affine-linear maps associated with it respectively. The set of functions $\Phi = \{ \phi_i \}_{i \in I}$ is called $regular$ $partition$ $of$ $unity$ $subordinate$ $to$ $\mathcal{Q}$,
  if
  %it satisfies the following:
  \begin{enumerate}
    \item $\phi_i \in C_c^\infty (O)$ with
          $\operatorname{supp} \phi_i \subset Q_i$ for all $i \in I$;

    \item $\sum_{i \in I} \phi_i \equiv 1$ on $O$; and

    \item $\sup_{i \in I}
               \| \, \partial^\alpha \phi_i^\natural \, \|_{L^\infty}
            < \infty $ for all $\alpha \in \mathbb{N}_0^2$, where
            $\phi^\natural :
            \mathbb{R}^d \to \mathbb{C},\xi\mapsto\phi_i(T_i \,\xi+b_i)$.
  \end{enumerate}
\end{defn}
 
\begin{defn}\label{moderate_weight}
The sequence $w = \{ w_i \}_{i \in I}$ of positive numbers is called $weight$. The weight is called $\mathcal{Q}$-$moderate$ where $\mathcal{Q} = \{ Q_i \}_{i \in I}$ stands for an almost structured cover of the open subset $O$ of $\mathbb{R}^d$, if there is such a positive number $C$ that $w_i \leqslant C \cdot w_{i'}$ for all such $i$ and $i' \in I$ that $Q_i \cap Q_{i'} \neq \emptyset$.
\end{defn} 
  
\begin{defn}\label{reservoir}
The topological dual $Z'$ of the space $Z := \mathcal{F} (C_c^\infty (\mathbb{R}^d)) \subset \mathcal{S} (\mathbb{R}^d)$ equipped with the unique topology that makes the Fourier transform
  $\mathcal{F} : C_c^\infty (\mathbb{R}^d) \to Z$ into a homeomorphism will be referred to as $reservoir$.
\end{defn}

Such a definition of the reservoir ensures that the decomposition space that we shall construct with its aid will be complete \cite{Bytchenkoff_2019}. 

\begin{defn}\label{decomposition_space}
Let $\mathcal{Q} = \{ Q_i \}_{i \in I}$, $\Phi = \{ \phi_i \}_{i \in I}$, $w = \{ w_i \}_{i \in I}$ and $Z'$ be an almost structured cover of the open subset $O$ of $\mathbb{R}^d$, a regular partition of unity on $O \subset \mathbb{R}^d$ subordinate to $\mathcal{Q}$, a $\mathcal{Q}$-moderate
  weight and the reservoir, respectively, and let $p$ and $q \in (0,\infty]$. The set
  \begin{equation}
    \mathcal{D} (\mathcal{Q}, L^p, \ell_w^q)
    := \left\{
         g \in Z' : \| g \|_{\mathcal{D} (\mathcal{Q},L^p,\ell_w^q)} < \infty
       \right\} \,
      \label{eq:DecompositionSpace}
  \end{equation}
  equipped with the quasi-norm
  \begin{equation}
    \| g \|_{\mathcal{D} (\mathcal{Q}, L^p, \ell_w^q)}
    := \left\|
          \left(
            w_i \cdot \| \mathcal{F}^{-1} (\varphi_i \cdot \widehat{g}) \|_{L^p}
          \right)_{i \in I}
       \right\|_{\ell^q} \in [0,\infty] \,
    \label{eq:DecompositionSpaceNorm}
  \end{equation} is called the $decomposition$ $space$.
\end{defn}
  
We shall now specify covers, regular partitions of unity subordinate to them and moderate weights for the anisotropic homo- and inhomogeneous Besov spaces.

\subsection{Construction of anisotropic homogeneous Besov spaces}\label{Construction of anisotropic homogeneous Besov spaces}

\begin{defn}\label{expansive_matrix}
The matrix $A$ whose elements are real numbers and whose spectrum $\sigma (A)$ is such that
  \begin{equation}  
  \min \limits_{\lambda \in \sigma (A)} \vert \lambda \vert > 1
  \label{eq:DilationMatrix}
  \end{equation} is called $expansive$.
\end{defn}

\begin{defn}\label{homo_cover}
Let $Q_0$ be such be a compact subset of $[-1, 1]^d \setminus \{ 0 \}$ and $A$ such a  $d \times d$ expansive matrix that 
\begin{equation*}
\bigcup\limits_{i \in \mathbb{Z}} Q_i = {\mathbb{R}}^d  \setminus \{ 0 \}
\end{equation*} where $Q_i := A^i Q_0$, then the set $\dot{Q}_B := \{Q_i\}_{i \in \mathbb{Z}}$ will be referred to as $anisotropic$ $homogeneous$ $cover$ of ${\mathbb{R}}^d  \setminus \{ 0 \}$.
\end{defn}

The proof of existence of such a set $\dot{Q}_B := \{Q_i\}_{i \in \mathbb{Z}}$ and its being indeed an almost structured cover of ${\mathbb{R}}^d  \setminus \{ 0 \}$ can be found in
Lemma 5.2 in \cite{Cheshmavar_2016}.

\begin{defn}\label{homo_partition_unity}
Let $\phi \in C^{\infty} ({\mathbb{R}}^d  \setminus \{ 0 \})$ and $\operatorname{supp} \phi \subset Q_0$ with $Q_0$ as in Definition~\ref{homo_cover}, then the set $\dot{\Phi}_B := \{\phi_i\}_{i \in \mathbb{Z}}$ where
  \begin{equation}  
  \phi_i (\xi) := \phi ( A^{-i} \xi)
  \label{eq:homo_partition_unity}
  \end{equation} as $\xi \in {\mathbb{R}}^d$ and $i \in \mathbb{Z}$ will be referred to as $anisotropic$ $homogeneous$ $partition$ $of$ $unity$.
\end{defn}

From Remark 2.3 in \cite{Cheshmavar_2016} can be inferred that $\dot{\Phi}_B$ is a regular partition of unity on ${\mathbb{R}}^d  \setminus \{ 0 \}$ subordinate to $\dot{Q}_B$.

\begin{defn}\label{homo_weight}
Let $s \in \mathbb{Z}$. Then the set $\dot{w}_B := \{w_i\}_{i \in \mathbb{Z}}$ where
  \begin{equation}  
  w_i := \vert \det A \vert^{si} 
  \label{eq:homo_weight}
  \end{equation} for any $i \in \mathbb{Z}$ will be referred to as $anisotropic$ $homogeneous$ $weight$.
\end{defn}

\begin{lem}\label{lemma_homogeneous_moderate_weight}
The weight $\dot{w}_B$ is $\dot{Q}_B$-moderate.
\end{lem}
\noindent $Proof$. In Lemma 5.2 of \cite{Cheshmavar_2016} it was proved that $Q_i \cap Q_{i'} = \emptyset$ where $Q_i$ and $Q_{i'} \in \dot{Q}_B$ if $\vert i-i'\vert$ exceeds a certain finite integer $\Delta i$. Therefore
\begin{equation*}
\frac{w_i}{w_{i'}} = \vert \det A \vert^{s(i-i')} \leqslant \vert \det A \vert^{\vert s \vert \Delta i}
\end{equation*} as $Q_i \cap Q_{i'} \neq \emptyset$ and $\vert \det A\vert^{\vert s \vert \Delta i}$ can be taken as the positive number $C$ mentioned in Definition \ref{moderate_weight} of the moderate weight. $\Box$

\begin{defn}\label{homo_Besov}
The decomposition space as defined by \eqref{eq:DecompositionSpace} with the almost structured cover $\dot{Q}_B$, partition of unity $\dot{\Phi}_B$ subordinate to $\dot{Q}_B$ and $\dot{Q}_B$-moderate weight $\dot{w}_B$ will be referred to as
$anisotropic$ $homogeneous$  $Besov$ $space$ and denoted by $\dot{B}_{p, q}^\alpha$.
\end{defn}

From this definition one can deduce that the fact that a function $g$ belongs to the space $\dot{B}_{p, q}^\alpha$ does not, generally speaking, imply that the composition of functions $g \circ R$ where $R$ stands for a rotation matrix will also belong to it, which justifies the attributive $anisotropic$ in the name of the space. Only if we choose to use a scalar expansive matrix $A$ for generating the space $\dot{B}_{p, q}^\alpha$, the latter will be isotropic. This remark, as it will transpire later, also applies to anisotropic inhomogeneous Besov spaces to whose definition we now turn.

\subsection{Construction of anisotropic inhomogeneous Besov spaces}\label{Construction of anisotropic heterogeneous Besov spaces}

\begin{defn}\label{hetero_cover}
Let $Q_0$ be such a subset of $[-1, 1]^d$, $Q_1$ such a compact subset of $[-1, 1]^d \setminus \{ 0 \}$ and $A$ be such a $d \times d$ expansive matrix that
\begin{equation*}
\bigcup\limits_{i \in \mathbb{N}_0} Q_i = {\mathbb{R}}^d
\end{equation*} where $Q_i := A^{i-1} Q_1$ as $i \in \mathbb{N}$, then the set $Q_B := \{Q_i\}_{i \in \mathbb{N}_0}$ will be referred to as $anisotropic$ $inhomogeneous$ $cover$ of ${\mathbb{R}}^d$.
\end{defn}

The proof of existence of such a set $Q_B := \{Q_i\}_{i \in \mathbb{N}_0}$ and its being indeed an almost structured cover of ${\mathbb{R}}^d$ can be found in
Lemma 5.2 in \cite{Cheshmavar_2016}.

\begin{defn}\label{hetero_partition_unity}
Let $\phi_0$ and $\phi_1 \in C^{\infty} ({\mathbb{R}}^d  \setminus \{ 0 \})$,  $\operatorname{supp} \phi_0 \subset Q_0$ and $\operatorname{supp} \phi_1 \subset Q_1$ with $Q_0$ and $Q_1$ as in Definition~\ref{hetero_cover}, then the set $\Phi_B := \{\phi_i\}_{i \in \mathbb{N}_0}$ where
  \begin{equation}  
  \phi_i (\xi) := \phi_1 ( A^{-(i-1)} \xi)
  \label{eq:hetero_partition_unity}
  \end{equation} as $\xi \in {\mathbb{R}}^d$ and $i \in \mathbb{N}$ will be referred to as $anisotropic$ $inhomogeneous$ $partition$ $of$ $unity$.
\end{defn}

From Remark 2.3 in \cite{Cheshmavar_2016} can be inferred that $\Phi_B$ is a regular partition of unity on ${\mathbb{R}}^d$ subordinate to $Q_B$.

\begin{defn}\label{hetero_weight}
Let $s \in \mathbb{Z}$. Then the set $\dot{w}_B := \{w_i\}_{i \in \mathbb{N}_0}$ where $w_0 := 1$ and
  \begin{equation}  
  w_i := \vert \det A \vert^{s (i-1)} 
  \label{eq:hetero_weight}
  \end{equation} for any $i \in \mathbb{N}$ will be referred to as $anisotropic$ $inhomogeneous$ $weight$.
\end{defn}

\begin{lem}\label{lemma_heterogeneous_moderate_weight}
The weight ${w}_B$ is $Q_B$-moderate.
\end{lem}

The proof of this lemma is almost identical to that of \ref{lemma_homogeneous_moderate_weight}.

\begin{defn}\label{hetero_Besov}
The decomposition space as defined by \eqref{eq:DecompositionSpace} with the almost structured cover $Q_B$, partition of unity $\Phi_B$ subordinate to $Q_B$ and $Q_B$-moderate weight $w_B$ will be referred to as
$anisotropic$ $inhomogeneous$ $Besov$ $space$ and denoted by $B_{p, q}^\alpha$.
\end{defn}

\section{Construction of Banach frames for and atomic decompositions of the decomposition spaces}\label{Construction of Banach frames for and atomic decompositions of the decomposition spaces}

\noindent The decomposition space is a example of the quasi-Banach space, i. e.
the complete quasi-normed vector space. Banach frames and sets of atoms provide the quasi-Banach space with those basic building blocks into which any of the element of the space can be decomposed or from which an element of the space can be synthesised. Here are the definitions of these two notions \cite{Groechenig_1991}, along with one axillary definition.

\begin{defn}\label{Solid_subspace}
Let $ A \subset {\mathbb{C}}^I$ where ${\mathbb{C}}^I$ stands for the quasi-Banach space consisting of sequences of complex numbers indexed by $i \in I$ and let $a = \{ a_i \}_{i \in I} \in A$. If $\vert a'_i \vert \leqslant \vert a_i \vert$ for all $i \in I$ implies that $a' = \{ a'_i \}_{i \in I} \in A$ and $\vert \vert a' \vert \vert _{{\mathbb{C}}^I} \leqslant \vert \vert a \vert \vert_{{\mathbb{C}}^I}$, then $A$ is called $solid$.
\end{defn}

\begin{defn}\label{Banach_frame}
A set $\{\psi_i\}_{i \in I}$ in the dual space $X'$ of a quasi-Banach space $X$ is called $Banach$ $frame$ for $X$ if there is a well-defined bounded map, called $analysis$ $operator$, $A : X \rightarrow x , f \mapsto \{\langle \psi_i , f \rangle\}_{i \in I}$ where $x := \{\{\langle \psi_i , f \rangle\}_{i \in I} : f \in X\}$ is a solid quasi-Banach subspace of ${\mathbb{C}}^I$ and there is such a bounded linear map $A_l^{-1} : x \rightarrow X$ that $A_l^{-1} \circ A = I_X$ where $I_X$ stands for an identity operator on $X$.
\end{defn}

\begin{defn}\label{atoms}
A set $\{\phi_i\}_{i \in I}$ in a quasi-Banach space $X$ is called $set$ $of$ $atoms$ in $X$ if there is a well-defined bounded map, called $synthesis$ $operator$, $S : x \rightarrow X , \{c_i\}_{i \in I} \mapsto  \sum_{i \in I} c_i \phi_i$ where the coefficient space $x := \{c_i\}_{i \in I}$ associated with $\{\phi_i\}_{i \in I}$ is a solid subspace of ${\mathbb{C}}^I$ and there is such a bounded linear map $S_r^{-1} : X \mapsto x$ that $S \circ S_r^{-1} = I_X$ where $I_X$ stands for an identity operator on $X$. The series expansion $g = \sum_{i \in I} c_i \phi_i$ of a given function $g \in X$ where $\{\phi_i\}_{i \in I}$ is a set of atoms is called $atomic$ $decomposition$ of $g$.
\end{defn}

To construct Banach frames and sets of atoms for the anisotropic Besov spaces $\dot{B}_{p, q}^\alpha$ and $B_{p, q}^\alpha$ we shall use the following two concepts and two theorems.

\begin{defn}\label{GSIS}
Let $\delta > 0$, $\{ T_i \bullet + b_i \}_{i \in I}$ be a denumerable set of invertible affine-linear maps on $\mathbb{R}^d$ and $\{ \phi_n \}_{n=1}^N$ a finite set of square integrable functions on $\mathbb{R}^d$. Then
  \begin{equation}  
  \Psi := \{ L_{\delta T_i^{-t} k} \psi_i (t) \}_{i \in I, \ k \in \mathbb{R}^d}
  \label{eq:Psi}
  \end{equation} where
  \begin{equation} 
  L_x f(z) := f(z-x) \,
    \quad \text{and} \quad \, 
  \psi_i (t) := \vert \det T_i \vert^{1/2} M_{b_i} \phi_{n_i} (T_i^t t)
  \label{eq:Psi_2}
  \end{equation} where $M_s f(z) := e^{2 \pi i s z} f(z)$ and $\phi_{n_i} \in \{ \phi_n \}_{n=1}^N$ is called $generalised$ $shift$-$invariant$  $system$.
\end{defn}

\begin{defn}\label{coefficient_space}
Let $p$ and $q \in (0, \infty ]$ and $w = \{ w_i \}_{i \in I}$ be a weight. Then
\begin{equation} 
  C_w^{p,q} := \left\{ (c_k^{(i)})_{i \in I, \ k \in \mathbb{Z}^d} \in \mathbb{C}^{I \times \mathbb{Z}^d}
 : 
         \| c \|_{C_w^{p,q}}
         := \left\|
              \left(
                | \det T_i \, |^{\frac{1}{2} - \frac{1}{p}}
                \cdot w_i
                \cdot \| (c_k^{(i)})_{k \in \mathbb{Z}^d} \|_{\ell^p}
              \right)_{i \in I}
            \right\|_{\ell^q}
         < \infty
       \right\}
\end{equation} is called $coefficient$ $space$ $associated$ $with$ $\Psi$.
\end{defn}

The functions $\{ \phi_n \}_{n=1}^N$ in Definition~\ref{GSIS} can be regarded as prototypes of all the functions in $\Psi$. The next two theorems specify the conditions on the prototypes under which $\Psi$ will constitute a Banach frame or a set of atoms for the decomposition space $\mathcal{D} (\mathcal{Q}, L^p, \ell_w^q)$ respectively.

\begin{thm}\label{theorem_Banach_frame}
Let $\epsilon$, $p_0$ and $q_0 \in (0,1]$, $p \in [p_0, \infty]$, $q \in [q_0, \infty]$ and $\Phi = \{ \phi_i \}_{i \in I}$ and $w = (w_i)_{i \in I}$  be respectively a regular partition of unity subordinate to $\mathcal{Q}$ and a $\mathcal{Q}$-moderate weight where $\mathcal{Q}$ stands for the almost structured cover of an open subset $O$ of $\mathbb{R}^d$
whose elements $\{ Q_i \}_{i \in I}$ are generated by having the invertible affine-linear transformations $\{ T_i \bullet + b_i \}_{i \in I}$ act on elements of the finite set $\{ Q'_n \}_{n=1}^N$ of non-empty open and bounded subsets $Q'_i$ of $\mathbb{R}^d$. Cf Definition~\ref{almost_struct_cover},~\ref{reg_partition_unity} and~\ref{moderate_weight}. Furthermore  let all the elements of the finite set $\{ \phi_n \}_{n=1}^N$ of square integrable functions on $\mathbb{R}^d$, introduced in Definition~\ref{GSIS}, satisfy the following conditions:
  \begin{enumerate}
    \item
      $\hat{\phi}_n \in C^\infty (\mathbb{R}^d)$;  
    \item
      $\hat{\phi}_n$ and all its partial derivatives are of polynomial growth at most;  
    \item
       $\hat{\phi}_n (\xi) \neq 0$ as $\xi \in \overline{Q'_n}$;
    \item \label{enu:BanachFrameSpaceSideDecay} $\phi_n \in C^1(\mathbb{R}^d)$
      and $\nabla \phi_n \in L^1(\mathbb{R}^d) \cap L^\infty(\mathbb{R}^d)$; and       
\item
\begin{equation}
    S_1 := \sup_{i \in I} \sum_{j \in I} N_{ij}^1 < \infty
    \qquad \text{and} \qquad
    S_2 := \sup_{j \in I} \sum_{i \in I} N_{ij}^1 < \infty \,  \\
\label{eq:S_1-S_2}
\end{equation} where
  
\begin{equation}
\begin{split}
        N_{ij}^1 &
        := \left( \frac{w_i}{w_j} \right)^\tau
           \cdot (1 + \| T_i^{-1} T_j \|)^\sigma \\
           & \hspace{2.5cm}
           \cdot \max_{|\beta| \leq 1}
                   \left(
                     |\det T_j|^{-1} \,
                     \int\limits_{Q_j}
                       \max_{|\alpha| \leq N}
                       \left|
                         \left[
                           \partial^\alpha
                             \widehat{\partial^\beta \phi_{n_i}}
                         \right] \big( T_i^{-1} (\xi - b_i) \big)
                       \right|
                     \, d \xi
                   \right)^\tau \, ,
\end{split}
\label{eq:M_Theorem_2}
\end{equation} $\tau := \min \{1,p,q\}$,
  \[
    \sigma := \tau \cdot \left( N + \frac{d}{\min \{1,p\}} \right) \,
    \quad \text{and} \quad
    N    := \left\lceil \frac{d+\epsilon}{\min\{1,p\}} \right\rceil.    
  \] 
  \end{enumerate}
\noindent Then there is such a 
          $C = C(\epsilon, p_0, q_0, d, \mathcal{Q}, w, \{ \phi_n \}_{n=1}^N) > 0$ that the generalised shift-invariant system  
  \begin{equation}  
  \widetilde{\Psi} := \{ \widetilde{\psi}_i (t) \}_{i \in I, \ k \in \mathbb{R}^d} := \{ L_{\delta T_i^{-t} k} \psi_i (-t) \}_{i \in I, \ k \in \mathbb{R}^d}
  \label{eq:Psi}
  \end{equation} where $\psi_i (t)$ as defined by~(\ref{eq:Psi_2}) with the coefficient space $C_w^{p,q}$ constitutes a Banach frame for the decomposition space $\mathcal{D} (\mathcal{Q}, L^p, \ell_w^q)$ as long as $\delta \in (0,\delta_0]$ where
 \[
            \delta_0
            = 1 \Big/ \Big[
                    1
                    + C
                      \cdot C_{Q,w}^4
                      \cdot \left( S_1^{1/\tau} + S_2^{1/\tau} \right)^2
                  \Big] \, ;
          \]
 and, in particular, 
  \begin{enumerate}
    \item the analysis operator
      \[
        A_{\widetilde{\Psi}} :
        \mathcal{D} (\mathcal{Q}, L^p, \ell_w^q) \to C_w^{p,q}, \
        f \mapsto \{
                    [\widetilde{\psi}_i \ast f] (\delta T_i^{-t} k)
                  \}_{i \in I, \ k \in \mathbb{Z}^d} \,\, ,
      \]
      where the convolution $\widetilde{\psi}_i \ast f$ is defined by
      \begin{equation}
        \big( \widetilde{\psi}_i \ast f \big)(t)
        = \sum_{j \in I}
            \mathcal{F}^{-1}
            \Big(
              \,\widehat{\widetilde{\psi}_i} \cdot \phi_j \cdot \widehat{f} \,
            \Big)(t) \,
        \label{eq:BanachFrameConvolutionDefinition}
      \end{equation}
      is well-defined and bounded as long as $\delta \in (0,1]$ and the series in \eqref{eq:BanachFrameConvolutionDefinition} converges normally in $L^\infty (\mathbb{R}^d)$. Moreover, if
      $f \in L^2(\mathbb{R}^d) \hookrightarrow S'(\mathbb{R}^d) \hookrightarrow Z'$,
      the convolution defined by
      \eqref{eq:BanachFrameConvolutionDefinition} agrees with its usual
      definition and
      \begin{equation}
        A_{\widetilde{\Psi}} f
        = \{
            \langle
               f, L_{\delta T_i^{-t}  k} \widetilde{\psi}_i
            \rangle
          \}_{i \in I, \ k\in \mathbb{Z}^d}
        \label{eq:AnalysisOperatorConsistency}
      \end{equation} for any $f \in L^2 (\mathbb{R}^d) \cap \mathcal{D} (\mathcal{Q}, L^p, \ell_w^q)$.
    \item
   there is such a map
         $A_{\widetilde{\Psi} l}^{-1} :
          C_w^{p,q} \to \mathcal{D} (\mathcal{Q}, L^p,\ell_w^q)$
         that $A_{\widetilde{\Psi} l}^{-1} \circ A_{\Psi}
         = id_{\mathcal{D} (\mathcal{Q}, L^p, \ell_w^q)}$ as long as $\delta \in (0, \delta_0]$.
  \end{enumerate}
\end{thm}

This is a slightly reformulated statement of Theorem 2.11 in \cite{Bytchenkoff_2019}. The proof of this theorem in its more general and simplified formulations can be found in \cite{Voigtlaender_2016} and \cite{Voigtlaender_2017} respectively.

\begin{thm}\label{theorem_atoms}
Let $\epsilon$, $p_0$ and $q_0 \in (0,1]$, $p \in [p_0, \infty]$, $q \in [q_0, \infty]$ and $w = \{ w_i \}_{i \in I}$ be $\mathcal{Q}$-moderate weight where $\mathcal{Q}$ stands for the almost structured cover of an open subset $O$ of $\mathbb{R}^d$
whose elements $\{ Q_i \}_{i \in I}$ are generated by having the invertible affine-linear transformations $\{ T_i \bullet + b_i \}_{i \in I}$ act on elements of the finite set $\{ Q'_n \}_{n=1}^N$ of non-empty open and bounded subsets $Q'_i$ of $\mathbb{R}^d$. Cf Definition~\ref{almost_struct_cover} and~\ref{moderate_weight}. Furthermore let all the elements of the finite set $\{ \phi_n \}_{n=1}^N$ of square integrable functions on $\mathbb{R}^d$, introduced in Definition~\ref{GSIS}, satisfy the following conditions:

  \begin{enumerate}
    \item
      $\hat{\phi}_n \in C^\infty (\mathbb{R}^d)$ where $\hat{\phi}_n$ stands for the Fourier transform of $\phi_n$;  
    \item
      $\hat{\phi}_n$ and all its partial derivatives are of polynomial growth at most;
    \item
       $\hat{\phi}_n (\xi) \neq 0$ as $\xi \in \overline{Q'_n}$;

    \item
       \[
         \sup_{t \in \mathbb{R}^d}
           \left[
              (1+|t|)^{\Lambda} \cdot | \phi_n (t) |
           \right]
         < \infty
       \] where $\Lambda := 1 + d / \min \{1,p\}$; and
    \item
there is such a set $\{ \rho_n \}_{n=1}^N$ of non-negative and absolutely integrable on $\mathbb{R}^d$ functions that 
  \begin{enumerate} 
   \item  
       \[
         \Big|
           \partial^\alpha \hat{\phi}_n (\xi)
         \Big|
         \leqslant \rho_n (\xi) \cdot (1+|\xi|)^{-(d+1+\epsilon)}
       \] as $\xi \in \mathbb{R}^d$ and for all such $\alpha \in \mathbb{N}_0^d$ that $|\alpha|\leqslant N$ with $N$ as defined in Theorem \ref{theorem_Banach_frame}; and
    \item
\begin{equation}
    S_3 := \sup_{i \in I} \sum_{j \in I} N_{ij}^2 < \infty
    \qquad \text{and} \qquad
    S_4 := \sup_{j \in I} \sum_{i \in I} N_{ij}^2 < \infty \,  \\
\label{eq:S_3-S_4}                     
\end{equation} where
\begin{equation}
\begin{split}
        N_{ij}^2 &
        := \left(
             \frac{w_i}{w_j} \cdot
             \left(
               |\det T_j| \,\big/\, |\det T_i|
             \right)^{\theta}
           \right)^{\!\tau}
           \!\! \cdot (1 + \| T_j^{-1} T_i \|)^{\sigma}  \\
           & \hspace{2.5cm}
           \cdot \left(
                   |\det T_i|^{-1}
                   \,\, \int\limits_{Q_i}
                            \rho_{n_j} \big(T_j^{-1}(\xi - b_j)\big)
                        \, d \xi
                 \right)^{\tau}
                 \hspace{0.25cm},
\end{split}                 
\label{eq:M_Theorem_1}                 
\end{equation}
% \[
%    \vartheta := (p^{-1} - 1)_+                                %\,\, , \qquad
%    %\qquad \text{and} \qquad
%    \tau := \min \{1,p,q\}                                     %\,
%  \]      
$\theta := (p^{-1} - 1)_+$, $\tau$  as defined in Theorem \ref{theorem_Banach_frame} and
    \[
    \sigma
    := \begin{cases}
         \tau \cdot (d+1)
         & \text{if } p \in [1,\infty] \\
         \tau \cdot \big(p^{-1}\cdot d + \lceil p^{-1}\cdot (d+\epsilon)\rceil\big)
         & \text{if } p \in (0,1)
       \end{cases}
       \hspace{0.25cm}.
  \]
\end{enumerate}
\end{enumerate}

\noindent Then there is such a 
          $C = C(\epsilon, p_0, q_0, d, \mathcal{Q}, w, \{ \phi_n \}_{n=1}^N) > 0$ that the generalised shift-invariant system $\Psi$ with the coefficient space $C_w^{p,q}$ constitutes a set of atoms for the decomposition space $\mathcal{D} (\mathcal{Q}, L^p, \ell_w^q)$ as long as $\delta \in (0,\delta_0]$ where
 \[
            \delta_0
            = \min \left\{
                      1,
                      \Big[
                        C \cdot \big(S_3^{1/\tau} + S_4^{1/\tau} \, \big)
                      \Big]^{-1}
                   \right\} \, ;
          \]
 and, in particular,  
  \begin{enumerate}
    \item the synthesis operator
      \begin{equation}
        S_{\Psi} :
        C_w^{p,q} \to \mathcal{D} (\mathcal{Q}, L^p, \ell_w^q), \
        \{c_k^{(i)}\}_{i \in I, k \in \mathbb{Z}^d}
        \mapsto \sum_{i \in I}
                  \sum_{k \in \mathbb{Z}^d}
                    \left[
                      c_k^{(i)}
                       L_{\delta T_i^{-t} k} \psi_i
                    \right]
\label{eq:synthesis}
\end{equation}
      is well-defined and bounded as long as $\delta \in (0,1]$,
      i.e. the sum over the index $k$ in~(\ref{eq:synthesis})
      converges absolutely for any index $i \in I$ to a function in
      $L_{S}^1 (\mathbb{R}^d) \cap S'(\mathbb{R}^d)$ and the sum of such functions over the index $i$ converges unconditionally in the
      weak$^\ast$ sense in $Z'$; and

   \item
   there is such a map
         $S_{\Psi r}^{-1} :
         \mathcal{D} (\mathcal{Q}, L^p,\ell_w^q) \to C_w^{p,q}$
         that $S_{\Psi} \circ S_{\Psi r}^{-1}
         = id_{(\mathcal{Q}, L^p, \ell_w^q)}$ as long as $\delta \in (0, \delta_0]$ and the action of $S_{\Psi l}^{-1}$ on any given
         $f \in \mathcal{D} (\mathcal{Q}, L^p, \ell_w^q)$ does not depend on $p$, $q$ and $w$.
  \end{enumerate} 
\end{thm}

This is a slightly reformulated statement of Theorem 2.10 in \cite{Bytchenkoff_2019}. The proof of this theorem in its more general and simplified formulations can be found in \cite{Voigtlaender_2016} and \cite{Voigtlaender_2017} respectively.
  
\section{Construction of Banach frames and atomic decompositions of the anisotropic Besov spaces}\label{Construction of Banach frames for and atomic decompositions of the anisotropic Besov spaces}

\subsection{Banach frames and atomic decompositions of anisotropic homogeneous Besov spaces}\label{Banach frames for and atomic decompositions of homogeneous Besov spaces}

\noindent First of all we define shift-invariant systems for the anisotropic homogeneous spaces and coefficient spaces associated with them. Here are these definitions.

\begin{defn}\label{homo_wavelets}
 Let $\delta > 0$, $I_d$ be $d \times d$ identity matrix, $A_i := A^i$, $\psi (t) \in L^1 (\mathbb{R}^d)$ and
  \begin{equation}  
  \psi_i (t) = \vert \det A_i \vert^{1/2} \psi (A_i^t t)
  \label{eq:psi_i}
  \end{equation} where $i \in \mathbb{Z}$, then the set
  \begin{equation}  
  \dot{\Psi}_B :=\{ L_{\delta A_i^{-t} k} \psi_i (t) \}_{i \in \mathbb{Z}, \ k \in \mathbb{Z}^d} = \{ \vert \det A_i \vert^{1/2} \psi (A_i^t t -\delta k) \}_{i \in \mathbb{Z}, \ k \in \mathbb{R}^d}
  \label{eq:AHBS}
  \end{equation} will be referred to as $anisotropic$ $homogeneous$ $Besov$ $wavelets$.
\end{defn}

\begin{defn}\label{homo_coefficient_space}
 Let $p$ and $q \in (0, \infty ]$, then
\begin{equation*} 
\dot{C}_{Bs}^{\hspace{0.2cm} p,q}  := \left\{ (c_k^{(i)})_{i \in \mathbb{Z}, \ k \in \mathbb{Z}^d} \in \mathbb{C}^{I \times \mathbb{Z}^d}
 : 
         \| c \|_{\dot{C}_{Bs}^{\hspace{0.2cm} p,q}}
         := \left\|
              \left(
                | \det A \, |^{i \left( \frac{1}{2} - \frac{1}{p} \right)}
                \cdot w_i
                \cdot \| (c_k^{(i)})_{k \in \mathbb{Z}^d} \|_{\ell^p}
              \right)_{i \in \mathbb{Z}}
            \right\|_{\ell^q}
         < \infty
       \right\}
\end{equation*} where $w_i \in \dot{w}_B$ will be referred to as $coefficient$ $space$ $associated$ $with$ $\dot{\Psi}_B$.
\end{defn}

%To construct Banach frames and sets of atoms for homogeneous Besov spaces we shall need the following lemma.

To prove the two theorems that establish the conditions on $\dot{\Psi}_B$ under which it will be a Banach frame or a set of atoms for $\dot{B}_{p, q}^\alpha$, we shall make use of an auxiliary lemma. Here is this lemma of ours as well as the lemma 2.2 from \cite{Bownik_2005} that we shall use to state and prove it.

\begin{lem}\label{lemma_Bownik}
Let $\lambda_-$ and $\lambda_+$ be such real numbers that
  \begin{equation}  
 1 < \lambda_- < \min \limits_{\lambda \in \sigma (A)} \vert \lambda \vert \leqslant \max \limits_{\lambda \in \sigma (A)} \vert \lambda \vert < \lambda_+
  \label{eq:DilationMatrix_2}
  \end{equation} where $\sigma (A)$ stands for the spectrum of a $d \times d$ expansive matrix $A$, then there is such a number $b>0$ that
\begin{equation}
   \frac{1}{b} \cdot \lambda_-^j \vert \xi \vert \leqslant \vert A^{j} \xi \vert \leqslant b \cdot \lambda_+^j \vert \xi \vert
   \label{eq:Ajxi}
\end{equation} and
\begin{equation}
   \frac{1}{b} \cdot \lambda_+^{-j} \vert \xi \vert \leqslant \vert A^{-j} \xi \vert \leqslant b \cdot \lambda_-^{-j} \vert \xi \vert
   \label{eq:A-jxi}
\end{equation} where $\xi \in {\mathbb{R}}^{d}$ and $j \in {\mathbb{N}}_0$.
\end{lem}

\begin{lem}\label{homo_lemma}
 Let $Q_0$ and $A$ be respectively such an open bounded subset of ${\mathbb{R}}^{d} \setminus \{ 0 \}$ and such a $d \times d$ expansive matrix that
\begin{equation}
   \bigcup\limits_{n \in \mathbb{Z}} Q_n = {\mathbb{R}}^{d}  \setminus \{ 0 \}
   \label{eq:cupQn_aniso}
\end{equation} where
\begin{equation}
{\lbrace Q_n := A^n Q_0 \rbrace}_{n \in \mathbb{Z}}
   \hspace{0.25cm}.
   \label{eq:Qn_aniso}
\end{equation} Furthermore let $a, \tau$ and $\sigma>0$,
\begin{equation}
L > \log_{\lambda_-} a  \hspace{0.25cm},
\label{eq:L}
\end{equation}
\begin{equation}
N > \log_{\lambda_-} \left( \frac{\lambda_+^{\sigma / \tau}}{a} \right)
\label{eq:N}
\end{equation} with $\lambda_-$ and $\lambda_+$ defined by~(\ref{eq:DilationMatrix_2}) and
\begin{equation}
   \vert \hat{\psi}  (\xi) \vert \leqslant C \cdot \min \lbrace 1, {\vert \xi \vert}^L \rbrace \left( 1+ \vert \xi \vert \right)^{-N}
   \label{eq:psi}
\end{equation} where $C>0$ and $\xi \in Q_0$. Then
\begin{equation}
   \sup_{n \in \mathbb{Z}} \sum_{m \in \mathbb{Z}} M_{mn} \leqslant S < \infty
   \label{eq:series21}
\end{equation}
and
\begin{equation}
   \sup_{m \in \mathbb{Z}} \sum_{n \in \mathbb{Z}}M_{mn} \leqslant S < \infty
   \label{eq:series22}
\end{equation} where
\begin{equation}
 M_{mn} :=
 a^{\tau (m-n)} \cdot \left( 1+ \Vert A^{n-m} \Vert \right)^\sigma \cdot \left[ \frac{1}{\vert Q_n \vert} \int\limits_{Q_n} \vert \hat{\psi} (A^{-m} \xi) \vert \, d\xi \right]^\tau
\label{eq:Mmn_}
\end{equation} and
\begin{equation}
S = C^\tau \cdot \left( 1+b \right)^\sigma \cdot \left(b \max \lbrace \frac{1}{r}, R \rbrace \right)^{\max \lbrace L, N \rbrace \tau} \cdot \left[ \frac{1}{1- \left( \frac{a}{\lambda_-^L } \right)^\tau} + \frac{2}{1- \frac{\lambda_+^\sigma}{\left( a \lambda_-^N \right)^\tau }} \right]
\label{eq:totalsum1}
\end{equation} with $b$ defined by~(\ref{eq:Ajxi}) and~(\ref{eq:A-jxi}) and where $r$ and $R$ are such that $0 < r \leqslant \vert \xi \vert \leqslant R < \infty$ as $\xi \in Q_0$.
\end{lem}

\noindent $Proof$. First of all we use~(\ref{eq:Qn_aniso}) to change the set over which the integration in~(\ref{eq:Mmn_}) is done from $Q_n$ to $Q_0$ and obtain
\begin{equation}
\begin{split}	
M_{mn} & = a^{\tau (m-n)} \cdot \left( 1+ \Vert A^{n-m} \Vert \right)^\sigma \cdot \left[ \frac{1}{\vert A^n Q_0 \vert} \int\limits_{Q_0} \vert \hat{\psi} (A^{-m} A^n \xi) \vert \vert A^n \vert \, d\xi \right]^\tau \\ 
& = a^{\tau (m-n)} \cdot \left( 1+ \Vert A^{n-m} \Vert \right)^\sigma \cdot \left[ \frac{1}{\vert Q_0 \vert} \int\limits_{Q_0} \vert \hat{\psi} (A^{n-m} \xi) \vert \, d\xi \right]^\tau
   \end{split}
   \label{eq:e0}
\end{equation} and then combine it with~(\ref{eq:psi}) to estimate $M_{mn}$ from above, namely
\begin{equation}
M_{mn} \leqslant a^{\tau (m-n)} \cdot \left( 1+ \Vert A^{n-m} \Vert \right)^\sigma \cdot \left[ \frac{C}{\vert Q_1 \vert} \int\limits_{Q_1} \frac{\min \lbrace 1, \vert A^{n-m} \xi \vert^L \rbrace}{\left( 1+ \vert A^{n-m} \xi \vert \right)^N} \, d\xi \right]^\tau \hspace{0.25cm}.
\label{eq:e}
\end{equation} Now we use~(\ref{eq:Ajxi}) and~(\ref{eq:A-jxi}) and distinguish two cases. If $n-m \geqslant 0$, then combining~(\ref{eq:e}) with~(\ref{eq:Ajxi}) results in
\begin{equation}
\begin{split}
M_{mn} & \leqslant a^{\tau (m-n)} \cdot \left( 1+ b \lambda_+^{n-m} \right)^\sigma \cdot \left[ \frac{C}{\vert Q_1 \vert} \int\limits_{Q_1} \frac{\min \lbrace 1, \left( b \lambda_+^{n-m} \vert \xi \vert \right)^L \rbrace}{\left( 1+ \frac{\lambda_-^{n-m} \vert \xi \vert}{b}  \right)^N} \, d\xi \right]^\tau \\
 &  \leqslant a^{\tau (m-n)} \cdot \left( 1+ b \lambda_+^{n-m} \right)^\sigma \cdot \left[ C \frac{\min \lbrace 1, \left( b \lambda_+^{n-m} R \right)^L \rbrace}{\left( 1+ \frac{\lambda_-^{n-m} r}{b}  \right)^N}  \right]^\tau
\hspace{0.25cm}.
\end{split}
\label{eq:c1}
\end{equation} Otherwise, if $n-m \leqslant 0$, then combining~(\ref{eq:e}) with~(\ref{eq:A-jxi}) results in
\begin{equation}
\begin{split}
M_{mn} & \leqslant a^{\tau (m-n)} \cdot \left( 1+ b \lambda_-^{n-m} \right)^\sigma \cdot \left[ \frac{C}{\vert Q_1 \vert} \int\limits_{Q_1} \frac{\min \lbrace 1, \left( b \lambda_-^{n-m} \vert \xi \vert \right)^L \rbrace}{\left( 1+ \frac{\lambda_+^{n-m} \vert \xi \vert}{b}  \right)^N} \, d\xi \right]^\tau \\
& \leqslant a^{\tau (m-n)} \cdot \left( 1+ b \lambda_-^{n-m} \right)^\sigma \cdot \left[ C \frac{\min \lbrace 1, \left( b \lambda_-^{n-m} R \right)^L \rbrace}{\left( 1+ \frac{\lambda_+^{n-m} r}{b}  \right)^N} \right]^\tau  \hspace{0.25cm}.
\end{split}
\label{eq:g}
\end{equation} 

To prove that the series~(\ref{eq:series21}) converges we shall divide it into several series and prove that each of them converges. In doing so we distinguish two further cases. If $n \geqslant 0$, then we divide the series~(\ref{eq:series21}) into two, namely
\begin{equation}
\sum_{m=-\infty}^{\infty} M_{mn} = \sum_{m=-\infty}^{n} M_{mn} + \sum_{m=n+1}^{\infty} M_{mn}
\hspace{0.25cm},
\label{eq:c_}
\end{equation} so that $n-m \geqslant 0$ and therefore $M_{mn}$ can be estimated as in~(\ref{eq:c1}) in the former series and $n-m \leqslant 0$ and therefore $M_{mn}$ can be estimated as in~(\ref{eq:g}) in the latter series. Using the substitution $m' := n-m$ in the former series in~(\ref{eq:c_}) and changing the order of summation results in
\begin{equation}
\sum_{m=-\infty}^{n} M_{mn} \leqslant \sum_{m'=0}^{\infty} M_{m'}
\label{eq:1000}
\end{equation} where
\begin{equation}
\begin{split}
M_{m'} & := \left( 1+b \right)^\sigma \cdot C^\tau \cdot \left( \frac{b}{r} \right)^{N \tau} \cdot \left[ \frac{\lambda_+^\sigma}{\left( a \lambda_-^N \right)^\tau } \right]^{m'} \geqslant a^{- \tau m'} \cdot \left( 1+ b \lambda_+^{m'} \right)^\sigma \cdot \left[ \frac{C}{\left( 1+ \frac{\lambda_-^{m'} r}{b}  \right)^N}  \right]^\tau \\ & \geqslant a^{- \tau m'} \cdot \left( 1+ b \lambda_+^{m'} \right)^\sigma \cdot \left[ C \frac{\min \lbrace 1, \left( b \lambda_+^{m'} R \right)^L \rbrace}{\left( 1+ \frac{\lambda_-^{m'} r}{b}  \right)^N}  \right]^\tau = M_{mn}
\end{split}
\hspace{0.05cm}.
\label{eq:Mm'_1}
\end{equation} The geometrical series on the right-hand side of~(\ref{eq:1000}) would converge, should its general term $M_{m'}$ satisfy the criterion
\begin{equation}
\lim_{m' \to \infty} \sup \frac{M_{m'+1}}{M_m'} < 1 
\label{eq:criterion}
\end{equation} or, in other words, 
\begin{equation}
\lim_{m' \to \infty} \sup a^{- \tau} \frac{\lambda_+^\sigma}{\lambda_-^{N \tau}} < 1 \hspace{0.25cm}.
\label{eq:criter1}
\end{equation} This holds if, as assumed in this lemma,~(\ref{eq:N}) does. Under this condition
\begin{equation}
\begin{array}{c}
\sum\limits_{m=-\infty}^{n} M_{mn} \leqslant \sum\limits_{m'=0}^{\infty} M_{m'} = \left( 1+b \right)^\sigma  \cdot C^\tau \cdot \left( \frac{b}{r} \right)^{N \tau} \cdot \frac{1}{1- \frac{\lambda_+^\sigma}{\left( a \lambda_-^N \right)^\tau }} \hspace{0.25cm}.
\end{array}
\label{eq:sum1}
\end{equation}

Using the substitution $m' := m-n$ in the latter series in~(\ref{eq:c_}) results in
\begin{equation}
\sum_{m=n+1}^{\infty} M_{mn} \leqslant \sum_{m'=1}^{\infty} M_{m'}  \leqslant \sum_{m'=0}^{\infty} M_{m'}
\label{eq:1001}
\end{equation} where
\begin{equation}
\begin{split}
M_{m'} & := \left( 1+b \right)^\sigma \cdot C^\tau \left( b R \right)^{L \tau} \cdot \left[ \left( \frac{a}{\lambda_-^L } \right)^\tau \cdot \right]^{m'} \geqslant C^\tau \cdot (1+b)^\sigma \cdot a^{\tau m'} \cdot (b \lambda_-^{-m'} R )^{L \tau} \\ & \geqslant
a^{\tau m'} \cdot \left( 1+ b \lambda_-^{-m'} \right)^\sigma \cdot \left[ C \frac{\min \lbrace 1, \left( b \lambda_-^{-m'} R \right)^L \rbrace}{\left( 1+ \frac{\lambda_+^{-m'} r}{b}  \right)^N} \right]^\tau = M_{mn}
 \hspace{0.25cm}.
 \end{split}
\label{eq:Mm'_2}
\end{equation} The geometrical series on the right-hand side of~(\ref{eq:1001}) would converge, should its general term $M_{m'}$ satisfy the criterion~(\ref{eq:criterion}), namely
\begin{equation}
\left( \frac{a}{\lambda_-^L} \right)^\tau < 1 \hspace{0.25cm}.
\end{equation} This holds if, as assumed in this lemma,~(\ref{eq:L}) does. Under this condition
\begin{equation}
\sum_{m=n+1}^{\infty} M_{mn} \leqslant \sum_{m'=0}^{\infty} M_{m'} = \frac{\left( 1+b \right)^\sigma C^\tau \left( b R \right)^{L \tau}}{1- \left( \frac{a}{\lambda_-^L } \right)^\tau} \hspace{0.25cm}.
\label{eq:sum2}
\end{equation}

If $n < 0$, then we divide the series~(\ref{eq:series21}) into three, namely
\begin{equation}
\sum_{m=-\infty}^{\infty} M_{mn} = \sum_{m=-\infty}^{n} M_{mn} + \sum_{m=n+1}^{-(n+1)} M_{mn}  + \sum_{m=-n}^{\infty} M_{mn}
\hspace{0.25cm},
\label{eq:c__}
\end{equation} so that $n-m \geqslant 0$ and therefore $M_{mn}$ can be estimated as in~(\ref{eq:c1}) in the first and third series and $n-m \leqslant 0$ and therefore $M_{mn}$ can be estimated as in~(\ref{eq:g}) in the second series. Using the substitution $m' := n-m$ in the first series in~(\ref{eq:c__}) and changing the order of summation results in the estimate identical to~(\ref{eq:1000}), where the geometrical series on the right-hand side will, as we already know, converge as long as~(\ref{eq:N}) holds. Under this condition the sum of the first series in~(\ref{eq:c__}) was estimated in~(\ref{eq:sum1}). Since $n \leqslant 0$ using the substitution $m' := m-n$ in the second series in~(\ref{eq:c__}) results in
\begin{equation}
\sum_{m=n+1}^{-(n+1)} M_{mn} \leqslant \sum_{m=n}^{-(n+1)} M_{mn} \leqslant \sum_{m'=0}^{-(2n+1)} M_{m'}  \leqslant \sum_{m'=0}^{\infty} M_{m'}
\label{eq:1002}
\end{equation} where the geometrical series on the right-hand side is identical to that in~(\ref{eq:1001}). Thus the second series in~(\ref{eq:c__}) will converge if~(\ref{eq:L}) holds. Under this condition the sum of the second series in~(\ref{eq:c__}) was estimated in~(\ref{eq:sum2}). Using the substitution $m' := n-m$ in the third series in~(\ref{eq:c__}) results in
\begin{equation}
\sum_{m=-n}^{\infty} M_{mn} = \sum_{m'=-2n}^{\infty} M_{m'} \leqslant \sum_{m'=0}^{\infty} M_{m'}
\label{eq:1003}
\end{equation} where the geometrical series on the right-hand side in identical to that in~(\ref{eq:1000}). Thus the third series in~(\ref{eq:c__}) will converge if~(\ref{eq:N}) holds. Under this condition its sum was already estimated in~(\ref{eq:sum1}). This completes the proof of convergence of the series~(\ref{eq:series21}) with any $n \in \mathbb{Z}$. Combining~(\ref{eq:c_}) and~(\ref{eq:c__}) with the estimates~(\ref{eq:sum1}) and~(\ref{eq:sum2}) results in~(\ref{eq:totalsum1}).

Changing the order of summation in the series~(\ref{eq:series22}) converts it into series~(\ref{eq:series21}). Therefore the series~(\ref{eq:series22}) converges if, as supposed in the lemma, $N$ and $L$ satisfy~(\ref{eq:N}) and~(\ref{eq:L}) respectively. The estimate of its sum is given by~(\ref{eq:totalsum1}). $\Box$

The next theorem establishes the conditions on $\dot{\Psi}_B$ under which it will be a Banach frame for $\dot{B}_{p, q}^\alpha$.

\begin{thm}\label{theorem_homo_Banach_frame}
Let $\epsilon$, $p_0$ and $q_0 \in (0,1]$.  Moreover let $\phi \in L^1(\mathbb{R}^d)$
satisfy the following conditions:

  \begin{enumerate}
    \item
      $\hat{\phi} \in C^\infty (\mathbb{R}^d)$;  
    \item
      $\hat{\phi}$ and all its partial derivatives are of polynomial growth at most;
   \item
       $\hat{\phi} (\xi) \neq 0$ as $\xi \in \overline{Q}_0$, where $Q_0$ as in Definition \ref{homo_cover};      
      
    \item
       $\phi \in C^1 (\mathbb{R}^d)$ and $\nabla \phi \in L^1(\mathbb{R}^d) \cap L^\infty(\mathbb{R}^d)$; and
  \item    
\begin{equation}
         \Big|
           \partial^\alpha \widehat{\partial^\beta \phi} (\xi)
         \Big|
         = \widehat{\gamma}_1 ( \xi ) \leqslant C \min \lbrace 1, {\vert \xi \vert}^{L_1} \rbrace (1+|\xi|)^{- N_1}
\label{eq:gamma_1}
\end{equation} where $C$ stands for a constant,
\begin{equation}
L_1 > s \log_{\lambda_-} \left( \vert \det A \vert \right) \hspace{0.25cm},
\label{eq:L_1}
\end{equation} and
\begin{equation}
N_1 > \log_{\lambda_-} \left( \frac{\lambda_+^{\sigma / \tau}}{\vert \det A \vert^s} \right)
\label{eq:N_1}
\end{equation} as $\xi \in \mathbb{R}^d$ and for all such $\alpha$ and $\beta \in \mathbb{N}_0^d$ that $|\alpha|\leqslant N$ and $\vert \beta \vert \leqslant 1$.
       
  \end{enumerate} 
\noindent Then there is such a 
          $\delta_0 = \delta_0(\epsilon, p_0, q_0, d, A, \phi) > 0$ that the anisotropic homogeneous Besov wavelets $\dot{\Psi}_B$ with the coefficient space $\dot{C}_{Bs}^{\hspace{0.2cm} p,q}$ constitutes a Banach frame for the anisotropic homogeneous  Besov space $\dot{B}_{p, q}^\alpha (A)$ as long as $\delta \in (0,\delta_0]$.
\end{thm}

\noindent $Proof.$ The four assumptions of this theorem are nothing but those of Theorem~\ref{theorem_Banach_frame} formulated for the generalised shift-invariant system $\dot{\Psi}_B$ with the coefficient space $\dot{C}_{Bs}^{\hspace{0.2cm} p,q}$ pertaining to the cover $\dot{Q}_B$ of the set $O :=\mathbb{R}^d \setminus \{ 0 \}$ and weight $\dot{w}_B$ that form the space $\dot{B}_{p, q}^\alpha$. Furthermore we recollect that $\{ T_i \}_{i \in \mathbb{Z}} = \{ A^i \}_{i \in \mathbb{Z}}$ and  $\{ b_i \}_{i \in \mathbb{Z}} = 0$ as $\Psi =\dot{\Psi}_B$ and that $\{ w_i \}_{i \in \mathbb{Z}} = \{ \vert \det A \vert^{is} \}_{i \in \mathbb{Z}}$ with $s \in\mathbb{Z}$ as $C_w^{p,q} = \dot{C}_{Bs}^{\hspace{0.2cm} p,q}$. Therefore $N_{ji}^1$ defined by~(\ref{eq:M_Theorem_2}) becomes
 \begin{equation}
     \begin{split}
        N_{ij}^1
        & = \left( \vert \det A \vert^s \right)^{\tau (i-j)}
           \cdot (1 + \| A^{j-i} \|)^\sigma
           \cdot \max_{|\beta| \leq 1}
                   \left(
                     \vert \det A \vert^{-j} \,
                     \int\limits_{Q_j}
                       \max_{|\alpha| \leq N}
                       \left|
                         \left[
                           \partial^\alpha
                             \widehat{\partial^\beta \phi}
                         \right] \big( A^{-i} \xi \big)
                       \right|
                     \, d \xi
              \right)^\tau \\
  & = \left( \vert \det A \vert^s \right)^{\tau (i-j)}
           \cdot (1 + \| A^{j-i} \|)^\sigma
           \cdot \max_{|\beta| \leq 1}
                   \left(
                     \frac{2^d}{\vert Q_j \vert} \,
                     \int\limits_{Q_j}
                       \max_{|\alpha| \leq N}
                       \left|
                         \left[
                           \partial^\alpha
                             \widehat{\partial^\beta \phi}
                         \right] \big( A^{-i} \xi \big)
                       \right|
                     \, d \xi
                   \right)^\tau          
\end{split} 
     \label{eq:M_Theorem_4}                  
\end{equation} where we also noted that $ \vert Q_j \vert = \vert A^j Q \vert = \vert \det A^j \vert \cdot \vert Q \vert = \vert \det A \vert^j \cdot 2^d$. The last expression in~(\ref{eq:M_Theorem_4}) clearly equates to $2^d \cdot M_{mn}$ with $M_{mn}$ defined by~(\ref{eq:Mmn_}) if $m=i$, $n=j$ and $a = \vert \det A \vert^s$ and as long as $L_1$ and $N_1$ in~(\ref{eq:gamma_1}) are not smaller than $L$ and $N$ in~(\ref{eq:psi}) respectively. According to Lemma~\ref{homo_lemma} the series~(\ref{eq:series21}) and~(\ref{eq:series22}) converge on the assumptions~(\ref{eq:L}) and~(\ref{eq:N}). Therefore the series~(\ref{eq:S_1-S_2}) converge on the assumptions~(\ref{eq:L_1}) and~(\ref{eq:N_1}). In other words the fifth assumption of Theorem~\ref{theorem_Banach_frame} follows from the fifth assumption of the present theorem. $\Box$

The next theorem establishes the conditions on $\dot{\Psi}_B$ under which it will be a set of atoms for $\dot{B}_{p, q}^\alpha$.

\begin{thm}\label{theorem_homo_atoms}
Let $\epsilon$, $p_0$ and $q_0 \in (0,1]$. Moreover let $\phi \in L^1(\mathbb{R}^d)$
satisfy the following conditions:

  \begin{enumerate}
    \item
      $\hat{\phi} \in C^\infty (\mathbb{R}^d)$;  
    \item
      $\hat{\phi}$ and all its partial derivatives are of polynomial growth at most;
    \item
       $\hat{\phi} (\xi) \neq 0$ as $\xi \in \overline{Q}_0$, where $Q_0$ as in Definition  \ref{homo_cover};

    \item
       \[
         \sup_{t \in \mathbb{R}^d}
           \left[
              (1+|t|)^{\Lambda} \cdot | \phi (t) |
           \right]
         < \infty
       \] where $\Lambda := 1 + d / {p_0}$; and
  \item    
       \[
         \Big|
           \partial^\alpha \hat{\phi} (\xi)
         \Big|
         \leqslant \rho (\xi) \cdot (1+|\xi|)^{-(d+1+\epsilon)}
       \] where
   \[
    \rho : \mathbb{R}^d \to (0, \infty), \ \xi \mapsto C \min \lbrace 1, {\vert \xi \vert}^{L_2} \rbrace (1+|\xi|)^{- N_2}
    \hspace{0.25cm},
  \] $C$ stands for a constant,
\begin{equation}
L_2 > (\theta-s) \log_{\lambda_-} \left( \vert \det A \vert \right)
\label{eq:L_2}
\end{equation} and
\begin{equation}
N_2 > \log_{\lambda_-} \left( \frac{\lambda_+^{\sigma / \tau}}{ \vert \det A \vert^{\theta-s}} \right)
\label{eq:N_2}
\end{equation} as $\xi \in \mathbb{R}^d$ and for all such $\alpha \in \mathbb{N}_0^d$ that $|\alpha|\leqslant N$ where
   \[
    N    := \left\lceil \frac{d+\epsilon}{p_0} \right\rceil
       \hspace{0.25cm}.
  \] 
       
  \end{enumerate} 
\noindent Then there is such a 
          $\delta_0 = \delta_0(\epsilon, p_0, q_0, d, A, \phi) > 0$ that the anisotropic homogeneous Besov wavelets $\dot{\Psi}_B$ with the coefficient space $\dot{C}_{Bs}^{\hspace{0.2cm} p,q}$ constitutes a set of atoms for the anisotropic homogeneous  Besov space $\dot{B}_{p, q}^\alpha (A)$ as long as $\delta \in (0,\delta_0]$.
\end{thm}

\noindent $Proof.$ The first four assumptions of this theorem are nothing but those of Theorem~\ref{theorem_atoms} formulated for the generalised shift-invariant system $\dot{\Psi}_B$ with the coefficient space $\dot{C}_{Bs}^{\hspace{0.2cm} p,q}$ pertaining to the cover $\dot{Q}_B$ of the set $O :=\mathbb{R}^d \setminus \{ 0 \}$ and weight $\dot{w}_B$ that form the space $\dot{B}_{p, q}^\alpha$.

Furthermore, given that $\{ T_i \}_{i \in \mathbb{Z}} = \{ A^i \}_{i \in \mathbb{Z}}$ and  $\{ b_i \}_{i \in \mathbb{Z}} = 0$ as $\Psi =\dot{\Psi}_B$ and that $\{ w_i \}_{i \in \mathbb{Z}} = \{ \vert \det A \vert^{is} \}_{i \in \mathbb{Z}}$ with $s \in\mathbb{Z}$ as $C_w^{p,q} = \dot{C}_{Bs}^{\hspace{0.2cm} p,q}$, $N_{ij}^2$ defined by~(\ref{eq:M_Theorem_1}) becomes       \begin{equation}
     \begin{split}
        N_{ij}^2
        & = \left(
             \vert \det A \vert^{s(i-j)} \cdot
              \vert \det A \vert^{\theta (j-i)}
           \right)^{\!\tau}
           \!\! \cdot (1 + \| A^{i-j} \|)^{\sigma}
           \cdot \left(
                  \vert \det A \vert^{-i}
                   \,\, \int\limits_{Q_i}
                            \rho \big(A^{-j} \xi \big)
                        \, d \xi
               \right)^{\tau} \\
                 & = \left( \vert \det A \vert^{\theta-s} \right)^{\tau (j-i)}
           \!\! \cdot (1 + \| A^{i-j} \|)^{\sigma}
           \cdot \left(
                   \frac{2^d}{\vert Q_i \vert}
                   \,\, \int\limits_{Q_i}
                            \rho \big(A^{-j} \xi \big)
                        \, d \xi
                 \right)^{\tau}
     \end{split}
     \hspace{0.25cm}.
     \label{eq:M_Theorem_3}
     \end{equation} The last expression in~(\ref{eq:M_Theorem_3}) equates to $2^d \cdot M_{mn}$ with $M_{mn}$ defined by~(\ref{eq:Mmn_}) if $m=-i$, $n=-j$ and $a = \vert \det A \vert^{\theta - s}$ and as long as $L_1$ and $N_1$ in~(\ref{eq:gamma_1}) are not smaller than $L$ and $N$ in~(\ref{eq:psi}) respectively. According to Lemma~\ref{homo_lemma} the series~(\ref{eq:series21}) and~(\ref{eq:series22}) converge on the assumptions~(\ref{eq:L}) and~(\ref{eq:N}). Therefore the series~(\ref{eq:S_3-S_4}) converge on the assumptions~(\ref{eq:L_2}) and~(\ref{eq:N_2}). In other words the fifth assumption of Theorem~\ref{theorem_atoms} follows from the fifth assumption of the present theorem. $\Box$

\subsection{Banach frames and atomic decompositions of anisotropic inhomogeneous Besov spaces}\label{Banach frames for and atomic decompositions of heterogeneous Besov spaces}

\noindent We now define shift-invariant systems for the anisotropic inhomogeneous spaces and coefficient spaces associated with them. Here are these definitions.

\begin{defn}\label{hetero_wavelets}
Let $\delta > 0$, $I_d$ be $d \times d$ identity matrix, $A_i := A^i$, $\psi'$ and $\psi (t) \in L^1 (\mathbb{R}^d)$ and
  \begin{equation}  
  \psi_i (t) = \vert \det A_{i-1} \vert^{1/2} \psi (A_{i-1}^t t)
  \label{eq:psi_i}
  \end{equation} where $i \in \mathbb{N}$, then the set
 \begin{equation}
\begin{split} 
\Psi_B & := \{ L_{\delta I_d k} \psi' (t) \}_{k \in  \mathbb{Z}^d} \cup \{ L_{\delta A_i^{-t} k} \psi_i (t) \}_{i \in \mathbb{N}_0, \ k \in  \mathbb{Z}^d} \\ & =  \{ \psi' (t - \delta I_d k) \}_{k \in  \mathbb{Z}^d} \cup \{ \vert \det A_i \vert^{1/2} \psi (A_i^t t -\delta k) \}_{i \in \mathbb{N}_0, \ k \in  \mathbb{R}^d}
\end{split}
  \label{eq:AHeBS}
  \end{equation} will be referred to as $anisotropic$ $inhomogeneous$ $Besov$ $wavelets$.
\end{defn} 

\begin{defn}\label{hetero_coefficient_space}
 Let $p$ and $q \in (0, \infty ]$, then
\begin{equation*} 
\begin{split}
C_{Bs}^{\hspace{0.2cm} p,q} :=
 & \left\{ (c_k)_{k \in \mathbb{Z}^d} \in \mathbb{C}^{\mathbb{Z}^d}
 : 
\| (c_k^{(i)})_{k \in \mathbb{Z}^d} \|_{\ell^p}
         < \infty
       \right\} \\ &
 \cup \left\{ (c_k^{(i)})_{i \in \mathbb{N}_0, \ k \in \mathbb{Z}^d} \in \mathbb{C}^{I \times \mathbb{Z}^d}
 : 
         \| c \|_{C_{Bs}^{\hspace{0.2cm} p,q}}
         := \left\|
              \left(
                | \det A \, |^{i \left( \frac{1}{2} - \frac{1}{p} \right)}
                \cdot w_i
                \cdot \| (c_k^{(i)})_{k \in \mathbb{Z}^d} \|_{\ell^p}
              \right)_{i \in \mathbb{N}_0}
            \right\|_{\ell^q}
         < \infty
       \right\}
\end{split}       
\end{equation*} where $w_i \in w_B$ will be referred to as $coefficient$ $space$ $associated$ $with$ $\Psi_B$ respectively.
\end{defn}

To prove the two theorems that establish the conditions on  $\Psi_B$ under which it will be a Banach frame or a set of atoms for  $B_{p, q}^\alpha$, we shall use the following lemma.

\begin{lem}\label{hetero_lemma}
 Let $Q_1$ be an open bounded subset of ${\mathbb{R}}^{d}$ that does not include its origin and the matrix $A \in {\mathbb{R}}^{d\times d}$ with eigenvalues ${\lbrace \lambda_i \in {\mathbb{C} : 1< \lambda_- < \vert \lambda_i \vert < \lambda_+} \rbrace}_{i=1}^d$ be such that the denumerable set
\begin{equation}
   {\lbrace Q_n := A^{n-1} Q_1 \rbrace}_{n \in \mathbb{N}}
   \label{eq:Qn_iso}
\end{equation} covers $ {\mathbb{R}}^{d} \setminus Q_0$ where $Q_0$ is an open bounded set that includes the origin of the ${\mathbb{R}}^{d}$. Furthermore let $a, \tau$, $\sigma>0$, $L$, $N$ and $\hat{\psi}$ be as defined in Lemma \ref{homo_lemma},
\begin{equation}
K \geqslant \log_{\lambda_-} \left( \frac{\lambda_+^{\sigma/ \tau}}{a} \right)
\label{eq:K}
\end{equation} and
\begin{equation}
   \vert \hat{\phi}  (\xi) \vert \leqslant C \left( 1+ \vert \xi \vert \right)^{-K}
   \label{eq:phi}
\end{equation} where $C>0$ and $\xi \in Q_1$. Then
\begin{equation}
   \sup_{n \in {\mathbb{N}}_0} \sum_{m \in {\mathbb{N}}_0} M_{mn} \leqslant S_{1} < \infty
   \label{eq:series1}
\end{equation} and

\begin{equation}
   \sup_{m \in {\mathbb{N}}_0} \sum_{n \in {\mathbb{N}}_0} M_{mn} \leqslant S_{2} < \infty
   \label{eq:series2}
\end{equation} where
\begin{equation}
 M_{mn} :=
  \begin{cases} 
   \begin{array}{l}
    2^\sigma \cdot \left[ \frac{1}{\vert Q_0 \vert} \int\limits_{Q_0} \vert \hat{\phi} (\xi) \vert \, d\xi \right]^\tau \hspace{0.25cm}, \hspace{0.25cm} (m=n=0) \\
    a^{\tau m} \cdot \left( 1+ \Vert A^{-(m-1)} \Vert \right)^\sigma \cdot \left[ \frac{1}{\vert Q_0 \vert} \int\limits_{Q_0} \vert \hat{\psi} ( A^{-(m-1)}\xi) \vert \, d\xi \right]^\tau \hspace{0.25cm}, \hspace{0.25cm} (m \neq 0, n=0)\\
    a^{- \tau n} \cdot \left( 1+ \Vert A^{n-1} \Vert \right)^\sigma \cdot \left[ \frac{1}{\vert Q_n \vert} \int\limits_{Q_n} \vert \hat{\phi} (\xi) \vert \, d\xi \right]^\tau \hspace{0.25cm}, \hspace{0.25cm} (m=0, n \neq 0)\\
   a^{\tau (m-n)} \cdot \left( 1+ \Vert A^{n-m} \Vert \right)^\sigma \cdot \left[ \frac{1}{\vert Q_n \vert} \int\limits_{Q_n} \vert \hat{\psi} (A^{-(m-1)} \xi) \vert \, d\xi \right]^\tau \hspace{0.25cm}, \hspace{0.25cm} (m \neq 0, n \neq 0)
   \end{array}
  \end{cases} \hspace{0.25cm},
\label{eq:Mmn}
\end{equation} and
\begin{equation}
 S_1 :=
  \begin{cases} 
   \begin{array}{l}
    2^\sigma \cdot C^\tau + \frac{(1+b)^\sigma (b \lambda_- R_0 )^{L \tau} C^\tau}{1- \left( \frac{a}{\lambda_-^L} \right)^\tau} \hspace{0.25cm}, \hspace{0.25cm} (n=0) \\
    
    \begin{split}
    C^\tau \cdot \left( 1+b \right)^\sigma & \cdot \left(1+ \frac{b \lambda_-}{\min \lbrace 1, r\rbrace} \right)^{\max \lbrace N, K \rbrace \tau} \\ & \cdot \left( 1+ b R \right)^{L \tau} \cdot \left[ \frac{1}{1- \frac{\lambda_+^\sigma}{(a \lambda_-^K )^\tau}} + \frac{1}{1- \frac{\lambda_+^\sigma}{\left( a \lambda_-^N \right)^\tau }} +  \frac{1}{1- \left( \frac{a}{\lambda_-^L } \right)^\tau} \right] \hspace{0.25cm}, \hspace{0.25cm} (n \neq 0)
        \end{split}
   \end{array}
  \end{cases} 
\label{eq:S1}
\end{equation} and
\begin{equation}
 S_2 :=
  \begin{cases} 
   \begin{array}{l}
    2^\sigma \cdot C^\tau + C^\tau \cdot \left( \frac{1+b}{\lambda_+} \right)^\sigma \cdot \left( \frac{b \lambda_-}{r} \right)^{K \tau} \cdot \frac{1}{1- \frac{\lambda_+^\sigma}{(a \lambda_-^K )^\tau}} \hspace{0.25cm}, \hspace{0.25cm} (m=0) \\
    
    \begin{split}    
    C^\tau \cdot (1+b)^\sigma & \cdot (b \lambda_- \max \lbrace R, R_0 \rbrace )^{L \tau}  \\ & \cdot \left( \frac{1+b}{\min \lbrace 1, r\rbrace } \right)^{N \tau} \cdot \left[ \frac{1}{1- \frac{\lambda_+^\sigma}{(a \lambda_-^N )^\tau }} + \frac{2}{1- \left( \frac{a}{\lambda_-^L} \right)^\tau}\right] \hspace{0.25cm}, \hspace{0.25cm} (m \neq 0)
        \end{split}    
   \end{array}
  \end{cases} 
\label{eq:S2}
\end{equation} with $b$ defined by~(\ref{eq:Ajxi}) and~(\ref{eq:A-jxi}) and where $R_0$, $r$ and $R$ are such that $0 \leqslant \vert \xi \vert < R_0$ as $\xi \in Q_0$ and $0 < r \leqslant \vert \xi \vert \leqslant R < \infty$ as $\xi \in Q_1$.
\end{lem}

\noindent \textit{Proof}. First of all we deal with the series~(\ref{eq:series1}) and consider the case where $n=0$. To do so we rewrite it as
\begin{equation}
\sum_{m=0}^{\infty} M_{mn} = M_{00} + \sum_{m=1}^{\infty} M_{m0}
 \hspace{0.25cm}.
\label{eq:ii}
\end{equation} Using~(\ref{eq:phi}) into~(\ref{eq:Mmn}) we estimate $M_{00}$ from above, namely
\begin{equation}
   M_{00} \leqslant 2^\sigma \cdot C^\tau \hspace{0.25cm}.
   \label{eq:M00}
\end{equation} Thus $M_{00}$ is bounded independently of $K$. Now we use~(\ref{eq:Ajxi}) and~(\ref{eq:A-jxi}) to estimate the general term of the series in~(\ref{eq:ii}) from above, namely
\begin{equation}
\begin{split}
  M_{m0} & \leqslant a^{\tau m} \cdot \left( 1+ b \lambda_-^{1-m} \right)^\sigma \cdot \left[ \frac{C}{\vert Q_0 \vert} \int\limits_{Q_0} \frac{\min \lbrace 1, \left( b \lambda_-^{1-m} \vert \xi \vert  \right)^L \rbrace}{\left( 1+ \frac{\lambda_+^{1-m} \vert \xi \vert }{b} \right)^N}\ d\xi \right]^\tau \\
  &  \leqslant a^{\tau m} \cdot \left( 1+ b \lambda_-^{1-m} \right)^\sigma \cdot \left[ \frac{C}{\vert Q_0 \vert} \int\limits_{Q_0} \min \lbrace 1, \left( b \lambda_-^{1-m} R_0  \right)^L \rbrace d\xi \right]^\tau \\
 &  \leqslant a^{\tau m} \cdot \left( 1+ b \lambda_-^{1-m} \right)^\sigma \cdot \left[ C \min \lbrace 1, \left( b \lambda_-^{1-m} R_0  \right)^L \rbrace \right]^\tau \\
& \leqslant a^{\tau m} \cdot \left( 1+ b \right)^\sigma \cdot \left( b \lambda_-^{1-m} R_0  \right)^{L\tau} \cdot C^\tau =: M_m
  \end{split}
  \hspace{0.25cm}.
   \label{eq:a}
\end{equation} Therefore the series in~(\ref{eq:series1}) would converge, should $M_m$ satisfy the criterion~(\ref{eq:criterion}), namely
\begin{equation}
\lim_{m \to \infty} \sup \frac{a^{\tau (m+1)} \cdot \left( 1+ b \right)^\sigma \cdot \left( b \lambda_-^{-m} R_0  \right)^{L\tau} \cdot C^\tau}{a^{\tau m} \cdot \left( 1+ b \right)^\sigma \cdot \left( b \lambda_-^{1-m} R_0  \right)^{L\tau} \cdot C^\tau} \leqslant \left( \frac{a}{\lambda_-^L} \right)^\tau < 1 \hspace{0.25cm}.
\label{eq:b}
\end{equation} This holds if, as assumed in this lemma,~(\ref{eq:L}) does. Under this condition
\begin{equation}
\sum_{m=1}^{\infty} M_{m0} \leqslant (1+b)^\sigma \cdot (b \lambda_- R_0 )^{L \tau} \cdot C^\tau \cdot \sum_{m=0}^{\infty} \left[ \left( \frac{a}{\lambda_-^L} \right)^\tau \right]^m = \frac{(1+b)^\sigma \cdot (b \lambda_- R_0 )^{L \tau} \cdot C^\tau}{1- \left( \frac{a}{\lambda_-^L} \right)^\tau} \hspace{0.25cm}.
\label{eq:2.1-2_}
\end{equation} Combining~(\ref{eq:M00}) and~(\ref{eq:2.1-2_}) results in~(\ref{eq:S1}) as $n=0$.

Now we investigate the series~(\ref{eq:series1}) as $n \neq 0$ and to do so divide it into three parts, namely
\begin{equation}
\sum_{m=0}^{\infty} M_{mn} = M_{0n} + \sum_{m=1}^{n} M_{mn} + \sum_{m=n+1}^{\infty} M_{mn}
\hspace{0.25cm},
\label{eq:c}
\end{equation} and deal with them separately. We use~(\ref{eq:Qn_iso}) to change the set over witch the integration in $M_{0n}$ is done from $Q_n$ to $Q_1$ and obtain $$ M_{0n} = a^{- \tau n} \left( 1+ \Vert A^{n-1} \Vert \right)^\sigma  \left[ \frac{1}{\vert Q_1 \vert} \int\limits_{Q_1} \vert \hat{\phi} ( A^{n-1} \xi) \vert \, d\xi \right]^\tau$$ and then combine it with~(\ref{eq:phi}) and~(\ref{eq:Ajxi}) to obtain
\begin{equation}
\begin{split}	
M_{0n} & \leqslant a^{- \tau n} \cdot \left( 1+ \Vert A^{n-1} \Vert \right)^\sigma \cdot \left[ \frac{C}{\vert Q_1 \vert} \int\limits_{Q_1}  \frac{\, d\xi}{\left( 1+ \vert A^{n-1} \xi \vert \right)^K}  \right]^\tau  \\ & \leqslant a^{- \tau n} \cdot \left( 1+ b \lambda_+^{n-1} \right)^\sigma \cdot \left[ \frac{C}{\vert Q_1 \vert} \int\limits_{Q_1}  \frac{\, d\xi}{\left( 1+ \frac{\lambda_-^{n-1} \vert \xi \vert}{b} \right)^K} \right]^\tau \\ & \leqslant
a^{- \tau n} \cdot \left( 1+ b \lambda_+^{n-1} \right)^\sigma \cdot \left[ \frac{C}{\left( 1+ \frac{\lambda_-^{n-1} r}{b} \right)^K}  \right]^\tau \leqslant a^{- \tau n} \cdot \left( 1+ b \lambda_+^{n-1} \right)^\sigma \cdot  \left[ \frac{C}{\left( \frac{\lambda_-^{n-1} r}{b} \right)^K}  \right]^\tau \\ 
& \leqslant \left( \frac{1+b}{\lambda_+} \right)^\sigma \cdot C^\tau \left( \frac{b \lambda_-}{r} \right)^{K \tau} \cdot \left[ \frac{\lambda_+^\sigma}{\left( a \lambda_-^K \right)^\tau } \right]^n
=: M_n
 \end{split}
   \hspace{0.25cm}.
\label{eq:Mn}
\end{equation} The $M_n$ function in~(\ref{eq:Mn}) is bounded for any finite $K$ and $n$. Therefore we only have to find the condition on $K$ under which $M_n$ stays bounded as $n$ grows unlimitedly. We shall achieve this by finding such a $K$ that 
\begin{equation}
\lim_{n \to \infty} M_n =0 \hspace{0.25cm}.
\label{eq:d}
\end{equation} This holds if, as assumed in this lemma,~(\ref{eq:K}) does. Under this condition
\begin{equation}
M_{0n} \leqslant \sum\limits_{n=1}^{\infty} M_{0n}  \leqslant \sum\limits_{n=1}^{\infty} M_{n} = \left( \frac{1+b}{\lambda_+} \right)^\sigma \cdot C^\tau \cdot \left( \frac{b \lambda_-}{r} \right)^{K \tau} \cdot \frac{1}{1- \frac{\lambda_+^\sigma}{(a \lambda_-^K )^\tau}}  \hspace{0.25cm}.
\label{eq:2.2-1_}
\end{equation}

 Now we deal with the rest of the series~(\ref{eq:c}). Changing the set over which the integration in $M_{mn}$ is done from $Q_n$ to $Q_1$ results in expression identical to~(\ref{eq:e0}). Combining it with~(\ref{eq:psi}) gives the estimate $M_{mn}$ from above identical to that in~(\ref{eq:e}). Using~(\ref{eq:Ajxi}), the general term in the second part of~(\ref{eq:c}), in which $n \geqslant m$, results in its estimate identical to that in~(\ref{eq:c1}). Making the substitution $m' := n-m$ and changing the order of summation leads to
\begin{equation}
\sum_{m=1}^{n} M_{mn} \leqslant
\sum_{m'=0}^{n-1} M_{m'}  \leqslant
\sum_{m'=0}^{\infty} M_{m'}
\label{eq:c2}
\end{equation} where $M_{m'}$ as defined in~(\ref{eq:Mm'_1}). Therefore the series in the second term of~(\ref{eq:c}) will converge if, as assumed in this lemma,~(\ref{eq:N}) holds. Its sum was already estimated in~(\ref{eq:sum1}). Using~(\ref{eq:A-jxi}), the general term in the third part of~(\ref{eq:c}), in which $n \leqslant m$, results in its estimate identical to that in~(\ref{eq:g}). Making the substitution $m' := m-n$ leads to
\begin{equation}
\sum\limits_{m=n+1}^{\infty} M_{mn} = \sum\limits_{m' =1}^{\infty} M_{m'}  \leqslant \sum\limits_{m' =0}^{\infty} M_{m'}
\label{eq:h}
\end{equation} where $M_{m'}$ as difined in~(\ref{eq:Mm'_2}). Therefore the series in the third term of~(\ref{eq:c}) will converge if, as assumed in this lemma,~(\ref{eq:L}) holds. Its sum was already estimated in~(\ref{eq:sum2}). Combining~(\ref{eq:2.2-1_}),~(\ref{eq:sum1}) and~(\ref{eq:sum2}) leads to~(\ref{eq:S1}) as the upper bound of the whole series in~(\ref{eq:c}) as $n \neq 0$.

Now we investigate the series~(\ref{eq:series2}). If $m=0$ we divide it into two parts, i.e.
\begin{equation}
\sum_{n=0}^{\infty} M_{0n} = M_{00} + \sum_{n=1}^{\infty} M_{0n}
\label{eq:i}
\end{equation} where $M_{00}$ and $M_{0n}$ were already estimated in~(\ref{eq:M00}) and~(\ref{eq:Mn})  respectively. Therefore the series~(\ref{eq:i}) will converge if, as assumed in this lemma,~(\ref{eq:K}) holds. Under this condition combining~(\ref{eq:M00}) and~(\ref{eq:Mn}) results in~(\ref{eq:S2}) as $m=0$.

If $m \neq 0$ we divide the series~(\ref{eq:series2}) into three parts
\begin{equation}
\sum_{n=0}^{\infty} M_{mn} = M_{m0} + \sum_{n=1}^{m} M_{mn} + \sum_{n=m+1}^{\infty} M_{mn}
\label{eq:j}
\end{equation} The term $M_{m0}$ was already estimated in~(\ref{eq:2.1-2_}) and will be bounded as long as the condition on $L$ in~(\ref{eq:L}) holds. The second term in~(\ref{eq:j}) was estimated in~(\ref{eq:g}) and can be converted into a finite sum over $m'$ whose general term is identical to that of the series in~(\ref{eq:h}) by using the substitution $m' := m-n$. The finite sum being smaller than the sum of the corresponding series, it leads us again to the condition on $L$ identical to that in~(\ref{eq:L}). Thus the upper bound of the second term in~(\ref{eq:j}) is given by~(\ref{eq:sum2}). Similarly the third part in~(\ref{eq:j}) can be converted into the series whose general term is identical to that in~(\ref{eq:Mm'_1}) by using the substitution $m' := n-m$. This leads to the condition on $N$ identical to that in~(\ref{eq:N}). The upper bound of the third part in~(\ref{eq:j}) is given by~(\ref{eq:sum1}). Combining~(\ref{eq:a}),~(\ref{eq:sum2}) and~(\ref{eq:sum1}) results in~(\ref{eq:S2}) as $m \neq 0$.  $\Box$

The next theorem establishes the conditions on $\Psi_B$ under which it will be a Banach frame for $B_{p, q}^\alpha$.

\begin{thm}\label{theorem_hetero_Banach_frame}
Let $\epsilon$, $p_0$ and $q_0 \in (0,1]$. Moreover let $\phi_1$ and $\phi_2 \in L^1(\mathbb{R}^d)$
satisfy the following conditions:

  \begin{enumerate}
    \item
      $\hat{\phi}_1$ and $\hat{\phi}_2 \in C^\infty (\mathbb{R}^d)$;  
    \item
      $\hat{\phi}_1$ and $\hat{\phi}_2$ and all its partial derivatives are of polynomial growth at most;
    \item
     $\hat{\phi}_1 (\xi) \neq 0$ as $\xi \in \overline{Q}_0$ and $\hat{\phi}_2 (\xi) \neq 0$ as $\xi \in \overline{Q}_1$ where $Q_0$ and $Q_1$ as in Definition \ref{hetero_cover};      
      
    \item
       $\phi \in C^1 (\mathbb{R}^d)$ and $\nabla \phi \in L^1(\mathbb{R}^d) \cap L^\infty(\mathbb{R}^d)$; and
  \item    
  \begin{equation}
         \Big|
           \partial^\alpha \widehat{\partial^\beta \phi_1} (\xi)
         \Big|
         \leqslant C (1+|\xi|)^{- K_4}
\label{eq:phi_1}
\end{equation} with $C$ standing for a constant and
\begin{equation}
K_4 > \log_{\lambda_-} \left( \frac{\lambda_+^{\sigma / \tau}}{\vert \det A \vert^{s-\theta}} \right)
\label{eq:K_4}
\end{equation} and
  \begin{equation}
         \Big|
           \partial^\alpha \widehat{\partial^\beta \phi_2} (\xi)
         \Big|
         \leqslant C \min \lbrace 1, {\vert \xi \vert}^{L_4} \rbrace (1+|\xi|)^{- N_4}
\label{eq:phi_2}         
\end{equation} with $C$ standing for a constant,
\begin{equation}
L_4 > (s-\theta) \log_{\lambda_-} \left( \vert \det A \vert \right)
\label{eq:L_4}
\end{equation} and
\begin{equation}
N_4 > \log_{\lambda_-} \left( \frac{\lambda_+^{\sigma / \tau}}{\vert \det A \vert^{s-\theta}} \right)
\label{eq:N_4}
\end{equation} as $\xi \in \mathbb{R}^d$ and $l \in \{ 1, 2 \}$ and for all such $\alpha$ and $\beta \in \mathbb{N}_0^d$ that $|\alpha|\leqslant N$ and $\vert \beta \vert \leqslant 1$.     
  \end{enumerate} 
\noindent Then there is such a 
          $\delta_0 = \delta_0(\epsilon, p_0, q_0, d, A, \phi) > 0$ that the anisotropic heterogeneous Besov wavelets $\Psi_B$ with the coefficient space $C_{Bs}^{\hspace{0.2cm} p,q}$ constitutes a Banach frame for the anisotropic heterogeneous  Besov space $B_{p, q}^\alpha (A)$ as long as $\delta \in (0,\delta_0]$.
\end{thm}

\noindent $Proof.$ The first four assumptions of this theorem are nothing but those of Theorem~\ref{theorem_Banach_frame} formulated for the generalised shift-invariant system $\Psi_B$ with the coefficient space $C_{Bs}^{\hspace{0.2cm} p,q}$ pertaining to the cover $Q_B$ of $\mathbb{R}^d$ and weight $w_B$ that form the space $B_{p, q}^\alpha$. Furthermore we recollect that $T_0 = I_d$, $\{ T_i \}_{i \in \mathbb{N}} = \{ A^{i-1} \}_{i \in \mathbb{N}}$ and  $\{ b_i \}_{i \in \mathbb{N}_0} = 0$ as $\Psi =\Psi_B$ and that $w_0 = 1$ and $\{ w_i \}_{i \in \mathbb{N}} = \{\vert \det A \vert^{(i-1) s} \}_{i \in \mathbb{N}}$ with $s \in \mathbb{Z}$ as $C_w^{p,q} = C_{Bs}^{\hspace{0.2cm} p,q}$. Therefore $N_{ij}^1$ defined by~(\ref{eq:M_Theorem_2}) becomes
 \begin{equation}
 N_{ij}^1 \leqslant
  \begin{cases} 
%   \begin{array}{l}
\begin{split}
    2^\sigma \cdot \max_{|\beta| \leqslant 1}
                   \left(
                     \frac{2^d}{\vert Q_0 \vert } \,
                     \int\limits_{Q_0}
                       \max_{|\alpha| \leqslant N}
                       \left|
                         \left[
                           \partial^\alpha
                             \widehat{\partial^\beta \phi_{1}}
                         \right] (  \xi )
                       \right|
                     \, d \xi
                   \right)^\tau \hspace{0.25cm}, \hspace{0.25cm} (i=j=0)
\end{split}                   
                    \\
                    \\
\begin{split}                    
     ( \vert \det A \vert^s )^{\tau (i-1)} & \cdot (1 + \| A^{-(i-1)} \|)^\sigma \\
    & \cdot \max_{|\beta| \leqslant 1}
                   \left(
                     \frac{2^d}{\vert Q_0 \vert } \,
                     \int\limits_{Q_0}
                       \max_{|\alpha| \leqslant N}
                       \left|
                         \left[
                           \partial^\alpha
                             \widehat{\partial^\beta \phi_{2}}
                         \right] \big( A^{-(i-1)} \xi \big)
                       \right|
                     \, d \xi
                   \right)^\tau \hspace{0.25cm}, \hspace{0.25cm} (i \neq 0, \ j=0)
\end{split}                   
                   \\
                   \\
\begin{split}
   (\vert \det A \vert^s )^{\tau (1-j)} &  \cdot (1 + \| A^{j-1} \|)^\sigma \\
    & \cdot \max_{|\beta| \leqslant 1}
                   \left(
                     \frac{2^d}{\vert Q_j \vert } \,
                     \int\limits_{Q_j}
                       \max_{|\alpha| \leqslant N}
                       \left|
                         \left[
                           \partial^\alpha
                             \widehat{\partial^\beta \phi_{1}}
                         \right] ( \xi )
                       \right|
                     \, d \xi
                   \right)^\tau \hspace{0.25cm}, \hspace{0.25cm} (i=0, \ j \neq 0)
\end{split}                   
                   \\
                   \\
\begin{split}
   (\vert \det A \vert^s )^{\tau (i-j)} & \cdot (1 + \| A^{j-i} \|)^\sigma \\
   & \cdot \max_{|\beta| \leqslant 1}
                   \left(
                     \frac{2^d}{\vert Q_j \vert } \,
                     \int\limits_{Q_j}
                       \max_{|\alpha| \leqslant N}
                       \left|
                         \left[
                           \partial^\alpha
                             \widehat{\partial^\beta \phi_{2}}
                         \right] \big( A^{-(i-1)} \xi \big)
                       \right|
                     \, d \xi
                   \right)^\tau \hspace{0.25cm}, \hspace{0.25cm} (i \neq 0, \ j \neq 0)
%   \end{array}
\end{split}
  \end{cases}
\label{eq:N_Banach_frame_hetero}
\end{equation} where we also noted that $\vert Q_0 \vert \leqslant 2^d$ and $ \vert Q_j \vert = \vert A^{j-1} Q_1 \vert = \vert \det A^{j-1} \vert \cdot \vert Q_1 \vert \leqslant \vert \det A \vert^{j-1} \cdot 2^d$ as $j \in \mathbb{N}$. The expressions on the right hand of~(\ref{eq:N_Banach_frame_hetero}) can be further estimated to conclude that, for any $i$ and $j \in \mathbb{N}_0$,
\begin{equation}
N_{ij}^1 \leqslant 2^d \cdot \vert \det A \vert^{\vert s \vert \tau} \cdot M_{mn}
\label{eq:N_Banach_frame_hetero_2}
\end{equation} with $M_{mn}$ defined by~(\ref{eq:Mmn}) if $m=i$, $n=j$ and $a = \vert \det A \vert^s$ and as long as $K_3$ in~(\ref{eq:rho_1}) and $L_3$ and $N_3$ in~(\ref{eq:rho_2}) are not smaller than $K$ in~(\ref{eq:phi}) and $L$ and $N$ in~(\ref{eq:psi}) respectively. According to Lemma~\ref{hetero_lemma} the series~(\ref{eq:series1}) and~(\ref{eq:series2}) converge on the assumption~(\ref{eq:psi}) about $\hat{\psi} (\xi)$ with $L$ defined by~(\ref{eq:L}), $N$ defined by~(\ref{eq:N}) and $K$ defined by~(\ref{eq:K}). Therefore the series~(\ref{eq:S_1-S_2}) converge on the assumption~(\ref{eq:rho_1}) about $\rho_1 (\xi)$ with $K_3$ defined by~(\ref{eq:K_3}) and the assumption~(\ref{eq:rho_2}) about $\rho_2 (\xi)$ with $L_3$ defined by~(\ref{eq:L_3}) and $N_3$ defined by~(\ref{eq:N_3}). In other words the fifth assumption of Theorem~\ref{theorem_Banach_frame} follows from the fifth assumption of the present theorem. $\Box$

The next theorem establishes the conditions on $\Psi_B$ under which it will be a set of atoms for $B_{p, q}^\alpha$.

\begin{thm}\label{theorem_hetero_atoms}
Let $\epsilon$, $p_0$ and $q_0 \in (0,1]$. Moreover let $\phi_1$ and $\phi_2 \in L^1(\mathbb{R}^d)$
satisfy the following conditions:

  \begin{enumerate}
    \item
      $\hat{\phi}_1$ and $\hat{\phi}_2 \in C^\infty (\mathbb{R}^d)$;  
    \item
      $\hat{\phi}_1$ and $\hat{\phi}_2$ and all its partial derivatives are of polynomial growth at most;
    \item
$\hat{\phi}_1 (\xi) \neq 0$ as $\xi \in \overline{Q}_0$ and $\hat{\phi}_2 (\xi) \neq 0$ as $\xi \in \overline{Q}_1$ where $Q_0$ and $Q_1$ as in Definition \ref{hetero_cover};

    \item for each $l \in \{ 1, 2 \}$
       \[
         \sup_{t \in \mathbb{R}^d}
           \left[
              (1+|t|)^{\Lambda} \cdot | \phi_l (t) |
           \right]
         < \infty
       \] where $\Lambda := 1 + d / {p_0}$; and
  \item    
       \[
         \Big|
           \partial^\alpha \hat{\phi}_l (\xi)
         \Big|
         \leqslant \rho_l (\xi) \cdot (1+|\xi|)^{-(d+1+\epsilon)}
       \] where
\begin{equation}
    \rho_1 : \mathbb{R}^d \to (0, \infty), \ \xi \mapsto C (1+|\xi|)^{- K_3}
\label{eq:rho_1}
\end{equation} with $C$ standing for a constant and
\begin{equation}
K_3 > \log_{\lambda_-} \left( \frac{\lambda_+^{\sigma / \tau}}{\vert \det A \vert^s} \right)
\label{eq:K_3}
\end{equation} and 
\begin{equation}
    \rho_2 : \mathbb{R}^d \to (0, \infty), \ \xi \mapsto C \min \lbrace 1, {\vert \xi \vert}^{L_3} \rbrace (1+|\xi|)^{- N_3}
\label{eq:rho_2}
\end{equation} with $C$ standing for a constant,
\begin{equation}
L_3 > s \log_{\lambda_-} \left( \vert \det A \vert \right)
\label{eq:L_3}
\end{equation} and
\begin{equation}
N_3 > \log_{\lambda_-} \left( \frac{\lambda_+^{\sigma / \tau}}{\vert \det A \vert^s} \right)
\label{eq:N_3}
\end{equation} as $\xi \in \mathbb{R}^d$ and for all such $\alpha \in \mathbb{N}_0^d$ that $|\alpha|\leqslant N$ where
   \[
    N    := \left\lceil \frac{d+\epsilon}{p_0} \right\rceil
       \hspace{0.25cm}.
  \] 
       
  \end{enumerate} 
\noindent Then there is such a 
          $\delta_0 = \delta_0(\epsilon, p_0, q_0, d, A, \phi) > 0$ that the anisotropic heterogeneous Besov wavelets $\Psi_B$ with the coefficient space $C_{Bs}^{\hspace{0.2cm} p,q}$ constitutes a set of atoms for the anisotropic heterogeneous  Besov space $B_{p, q}^\alpha (A)$ as long as $\delta \in (0,\delta_0]$.
\end{thm}

\noindent $Proof.$ The first four assumptions of this theorem are nothing but those of Theorem~\ref{theorem_atoms} formulated for the generalised shift-invariant system $\Psi_B$ with the coefficient space $C_{Bs}^{\hspace{0.2cm} p,q}$ pertaining to the cover $Q_B$ of $\mathbb{R}^d$ and weight $w_B$ that form the space $B_{p, q}^\alpha$. Furthermore we recollect that $T_0 = I_d$, $\{ T_i \}_{i \in \mathbb{N}} = \{ A^{i-1} \}_{i \in \mathbb{N}}$ and  $\{ b_i \}_{i \in \mathbb{N}_0} = 0$ as $\Psi =\Psi_B$ and that $w_0 = 1$ and $\{ w_i \}_{i \in \mathbb{N}} = \{\vert \det A \vert^{(i-1) s} \}_{i \in \mathbb{N}}$ with $s \in \mathbb{Z}$ as $C_w^{p,q} = C_{Bs}^{\hspace{0.2cm} p,q}$. Therefore $N_{ij}^2$ defined by~(\ref{eq:M_Theorem_1}) becomes
\begin{equation}
 N_{ij}^2 \leqslant
  \begin{cases} 
\begin{split}
    2^\sigma \cdot \left(
                   \frac{2^d}{\vert Q_0 \vert }
                   \,\, \int\limits_{Q_0}
                            \rho_1 ( \xi )
                        \, d \xi
                 \right)^{\tau}\hspace{0.25cm}, \hspace{0.25cm} (i=j=0)
\end{split}                   
                    \\
\begin{split}                    
    (\vert \det A \vert^{s - \theta} )^{\tau (i-1)} \cdot (1 + \| A^{i-1} \|)^\sigma  \cdot \left(
                   \frac{2^d}{\vert Q_i \vert }
                   \,\, \int\limits_{Q_i}
                            \rho_1 ( \xi )
                        \, d \xi
                 \right)^{\tau} \hspace{0.25cm}, \hspace{0.25cm} (i \neq 0, \ j=0)
\end{split}                   
                   \\
\begin{split}
    ( \vert \det A \vert^{s - \theta} )^{\tau (1-j)}  \cdot (1 + \| A^{-(j-1)} \|)^\sigma \cdot \left(
                   \frac{2^d}{\vert Q_0 \vert }
                   \,\, \int\limits_{Q_0}
                            \rho_2 \big( A^{-(j-1)} \xi \big)
                        \, d \xi
                 \right)^{\tau} \hspace{0.25cm}, \hspace{0.25cm} (i=0, \ j \neq 0)
\end{split}                   
                   \\
\begin{split}
   (\vert \det A \vert^{s - \theta} )^{\tau (i-j)} \cdot (1 + \| A^{i-j} \|)^\sigma \cdot \left(
                   \frac{2^d}{\vert Q_i \vert }
                   \,\, \int\limits_{Q_i}
                            \rho_2 \big( A^{-(j-1)} \xi \big)
                        \, d \xi
                 \right)^{\tau} \hspace{0.25cm}, \hspace{0.25cm} (i \neq 0, \ j \neq 0)
\end{split}
  \end{cases}
\label{eq:N_atoms_hetero}
\end{equation}  where we also noted that $\vert Q_0 \vert \leqslant 2^d$ and $ \vert Q_j \vert = \vert A^{j-1} Q_1 \vert = \vert \det A^{j-1} \vert \cdot \vert Q_1 \vert \leqslant \vert \det A \vert^{j-1} \cdot 2^d$ as $j \in \mathbb{N}$. The expressions on the right hand of~(\ref{eq:N_atoms_hetero}) can be further estimated to conclude that, for any $i$ and $j \in \mathbb{N}_0$,
\begin{equation}
N_{ij}^1 \leqslant 2^d \cdot \vert \det A \vert^{\vert s - \theta \vert \tau} \cdot M_{mn}
\label{eq:N_Banach_frame_hetero_2}
\end{equation} with $M_{mn}$ defined by~(\ref{eq:Mmn}) if $m=j$, $n=i$ and $a = \vert \det A \vert^{s - \theta}$ and as long as $K_4$ in~(\ref{eq:phi_1}) and $L_4$ and $N_4$ in~(\ref{eq:phi_2}) are not smaller than $K$ in~(\ref{eq:phi}) and $L$ and $N$ in~(\ref{eq:psi}) respectively. According to Lemma~\ref{hetero_lemma} the series~(\ref{eq:series1}) and~(\ref{eq:series2}) converge on the assumption~(\ref{eq:psi}) about $\hat{\psi} (\xi)$ with $L$ defined by~(\ref{eq:L}), $N$ defined by~(\ref{eq:N}) and $K$ defined by~(\ref{eq:K}). Therefore the series~(\ref{eq:S_1-S_2}) converge on the assumption~(\ref{eq:phi_1}) about $\phi_1 (\xi)$ with $K_4$ defined by~(\ref{eq:K_4}) and the assumption~(\ref{eq:phi_2}) about $\phi_2 (\xi)$ with $L_4$ defined by~(\ref{eq:L_4}) and $N_4$ defined by~(\ref{eq:N_4}). In other words the fifth assumption of Theorem~\ref{theorem_atoms} follows from the fifth assumption of the present theorem. $\Box$

%We shall premise the proof of this theorem on Theorem 1 and Lemma 1.          

\section*{Acknowledgements}
\noindent I thank the Centre National de la Recherche Scientifique of France
and the Deutscher Akademischer Austauschdienst of Germany for their funding and Professor Gitta Kutyniok for her support of this work.


\section*{References}
\begin{thebibliography}{99}

\bibitem{Besov_1961}
  O. V. Besov,
  \textit{Study of one family of spaces of functions in connexion with theorems of embedding and extension}, Trudi MIAN SSSR (1961) 42-81.

\bibitem{Besov_1975}
  O. V. Besov, V. P. Ilin and S. M. Nikolski,
  \textit{Integral representations of functions and embedding theorems}, Nauka, Moscow (1975).

\bibitem{Borup_2007}
 L. Borup and M. Nielsen,
  \textit{Frame decomposition of decomposition spaces III},
  J. Fourier Anal. Appl. 13 (2007) 39–70.   

\bibitem{Bownik_2005}
  M. Bownik,
  \textit{Atomic and molecular decompositions of anisotropic Besov spaces},
  Math. Z. 250 (2005) 539-571.    

\bibitem{Bytchenkoff_2019}
  D. Bytchenkoff and F. Voigtlaender,
  \textit{Design and properties of wave packet smoothness spaces},
  J. Math. Pures Appl. 133 (2020) 185-262.

\bibitem{Cheshmavar_2016}
  J. Cheshmavar and H. Führ,
  \textit{A classification of anisotropic Besov spaces},
  Appl. Harmon. Comput. Anal. (in press) doi.org/10.1016/j.acha.2019.04.006.          
  
\bibitem{Daubechies_1992}
  I. Daubechies,
  \textit{Ten lectures on wavelets},
  Society for industrial and applied mathematics, Philadelphia, PA
(1992).            
  
\bibitem{Feichtinger_1985}
  H. G. Feichtinger and P. Gröbner,
  \textit{Banach spaces of distributions defined by decomposition method},
  Nachr. 123 (1985) 97-120.

\bibitem{Frazier_1985}
 M. Frazier, B. Jawerth
  \textit{Decomposition of Besov spaces},
  Indiana Univ. Math. J. 34 (1985) 777-799.
  
\bibitem{Frazier_1989}
 M. Frazier, B. Jawerth
  \textit{A discrete transform and decomposition of distribution spaces},
  J. Funct. Anal. 93 (1989) 34-170.
  
\bibitem{Groechenig_1991}
  K. Gröchenig,
  \textit{Describing functions: atomic decompositions versus frames},
  Mh. Math. 112 (1991) 1-41.
  
\bibitem{Nikolski_1969}
  S. M. Nikolski,
  \textit{Approximation of functions of several variables and embedding theorems}, Nauka, Moscow (1969).  
  
\bibitem{Sobolev_1950}
  S. L. Sobolev,
  \textit{Some applications of functional analysis in mathematical physics},
  Leningrad State University, Leningrad (1950).    
  
\bibitem{Triebel_2006}
  H. Triebel,
  \textit{Theory of function spaces III},
  Birkhäuser Verlag, Basel (2006). 
  
\bibitem{Voigtlaender_2016}
  F. Voigtlaender,
  \textit{Structured, compactly supported Banach frame decompositions of decomposition spaces},
  http://arxiv.org/abs/1612.08772. 
  
\bibitem{Voigtlaender_2017}
  F. Voigtlaender and A. Pein,
  \textit{Analysis vs. synthesis sparsity for $\alpha$-shealets},
  http://arxiv.org/abs/1702.03559.     
  
   
\end{thebibliography}
\end{document}